\documentclass[12pt,oneside,reqno]{amsart}
\usepackage{txfonts}
\usepackage{bbm}
\usepackage{amsmath}
\usepackage{graphicx}
\usepackage{mathrsfs}
\usepackage{stmaryrd}
\usepackage{amsfonts}
\usepackage{xcolor}
\usepackage{verbatim}
\usepackage{enumerate,amsmath,amssymb,amsthm}
\usepackage{hyperref}

\oddsidemargin=0mm \topmargin=-12mm \numberwithin{equation}{section}

\newcommand{\be}{\begin{eqnarray}}
\newcommand{\ee}{\end{eqnarray}}
\newcommand{\ce}{\begin{eqnarray*}}
\newcommand{\de}{\end{eqnarray*}}
\newtheorem{theorem}{Theorem}[section]
\newtheorem{lemma}[theorem]{Lemma}
\newtheorem{remark}[theorem]{Remark}
\newtheorem{definition}[theorem]{Definition}
\newtheorem{proposition}[theorem]{Proposition}
\newtheorem{Examples}[theorem]{Example}
\newtheorem{corollary}[theorem]{Corollary}

\def\bu{{\mathbf{u}}}
\def\e{{\mathrm{e}}}
\def\eps{\varepsilon}

\def\a{\alpha}

\def\p{\partial}

\def\[{{\Big[}}
\def\]{{\Big]}}
\def\<{{\langle}}
\def\>{{\rangle}}
\def\({{\Big(}}
\def\){{\Big)}}

\def\bx{{\mathbf{x}}}

\def\dif{{\mathord{{\rm d}}}}

\def\min{{\mathord{{\rm min}}}}

\def\no{\nonumber}
\def\={&\!\!=\!\!&}
\def\bt{\begin{theorem}}
\def\et{\end{theorem}}
\def\bl{\begin{lemma}}
\def\el{\end{lemma}}
\def\br{\begin{remark}}
\def\er{\end{remark}}

\def\bd{\begin{definition}}
\def\ed{\end{definition}}
\def\bp{\begin{proposition}}
\def\ep{\end{proposition}}
\def\bc{\begin{corollary}}
\def\ec{\end{corollary}}
\def\bx{\begin{Examples}}
\def\ex{\end{Examples}}

\def\cB{{\mathcal B}}
\def\cC{{\mathcal C}}
\def\cD{{\mathcal D}}

\def\cH{{\mathcal H}}

\def\mC{{\mathbb C}}

\def\mE{{\mathbb E}}

\def\mH{{\mathbb H}}
\def\mI{{\mathbb I}}

\def\mL{{\mathbb L}}

\def\mN{{\mathbb N}}

\def\mP{{\mathbb P}}

\def\mR{{\mathbb R}}

\def\bP{{\mathbf P}}
\def\bE{{\mathbf E}}

\def\sD{{\mathscr D}}

\def\sF{{\mathscr F}}

\def\sI{{\mathscr I}}

\def\sK{{\mathscr K}}
\def\sL{{\mathscr L}}
\def\sM{{\mathscr M}}

\def\sP{{\mathscr P}}

\def\ba{{\begin{align}}
\def\ea{\end{align}}}

\def\geq{\geqslant}
\def\leq{\leqslant}

\def\bP{{\mathbf P}}

\allowdisplaybreaks

\def\bt{\begin{theorem}}
\def\et{\end{theorem}}
\def\bl{\begin{lemma}}
\def\el{\end{lemma}}
\def\br{\begin{remark}}
\def\er{\end{remark}}
\def\bx{\begin{Examples}}
\def\ex{\end{Examples}}

\def\ba{\begin{align}}
\def\ea{\end{align}}

\def\bd{\begin{definition}}
\def\ed{\end{definition}}
\def\bp{\begin{proposition}}
\def\ep{\end{proposition}}
\def\bc{\begin{corollary}}
\def\ec{\end{corollary}}
\def\bpf{\begin{proof}}
\def\epf{\end{proof}}

\begin{document}

\title{Heat kernel and ergodicity of SDEs with distributional drifts}
\date{}
\author{Xicheng Zhang and Guohuan Zhao}

\address{Xicheng Zhang:
School of Mathematics and Statistics, Wuhan University,
Wuhan, Hubei 430072, P.R.China\\
Email: XichengZhang@gmail.com
 }

\address{Guohuan Zhao:
Applied Mathematics, Chinese Academy of Science,
Beijing, 100081, P.R.China\\
Email: zhaoguohuan@gmail.com
 }

\thanks{
Research of X. Zhang is partially supported by NNSFC grant of China (No. 11731009). 
Research of G. Zhao is partially supported by National Postdoctoral Program for Innovative Talents (201600182) of China.}

\begin{abstract}
In this paper we consider the following SDE with distributional drift $b$:
$$
\dif X_t=\sigma(X_t)\dif B_t+b(X_t)\dif t,\ X_0=x\in\mR^d,
$$
where $\sigma$ is a bounded continuous and uniformly non-degenerate 
$d\times d$-matrix-valued function and 
$B$ is a $d$-dimensional standard Brownian motion.
Let $\alpha\in(0,\frac{1}{2}]$, $p\in(\frac{d}{1-\alpha},\infty)$ and $\beta\in[\alpha,1]$, $q\in(\frac{d}{\beta},\infty)$.
Assume
 $\|(\mI-\Delta)^{-\alpha/2}b\|_p+\|(-\Delta)^{\beta/2}\sigma\|_q<\infty$. 
We show the existence and uniqueness of martingale solutions to the above SDE, 
and obtain sharp two-sided and gradient estimates of the heat kernel associated to the above SDE.
Moreover, we study the ergodicity and  global regularity of the invariant measures of the associated semigroup under some dissipative assumptions.


\end{abstract}

\maketitle \rm

\section{Introduction}

Let $\sD$ be the space of all smooth functions on $\mR^d$ with compact supports,
and let $\sD'$ be the dual space of $\sD$, which is also called distributional function space.
Consider the following stochastic differential equation (abbreviated as SDE) in $\mR^d$ with distributional drift $b\in\sD'$:
\begin{align}\label{SDE}
\dif X_t=\sigma(X_t)\dif B_t+b(X_t)\dif t,\ \ X_0=x\in \mR^d,
\end{align}
where $B$ is a $d$-dimensional standard Brownian motion on some complete filtered probability space $(\Omega,\sF,(\sF_t)_{t\geq 0},\bP)$, 
 and $\sigma:\mR^d\to\mR^d\otimes\mR^d$ 
is a bounded continuous and non-degenerate $d\times d$-matrix-valued function. 
Since $b$ may be not a real function, the drift term $b(X_t)\dif t$ in \eqref{SDE} 
does not make any sense in general. A quite natural definition of the solution to SDE \eqref{SDE} is that $X$ is a continuous $\sF_t$-adapted process 
and satisfies
\begin{align}\label{SDE00}
X_t=x+\int^t_0\sigma(X_s)\dif B_s+A^b_t \mbox{ with }A^b_t:=\lim_{n\to\infty}\int^t_0 b_n(X_s)\dif s,
\end{align}
where $(b_n)_{n\in\mN}$ is any mollifying approximation sequence of $b$, and the limit is taken in the sense of u.c.p 
(uniformly on compact subsets of time variable in probability).
Suppose now that $b\in H^{-\alpha,p}$ for some $\alpha>0$ and $p>1$, where $H^{-\alpha,p}$ is the usual Bessel potential space
(see Definition \ref{def-Bessel} below).
To show the existence of the above limit, 
one possible way 
is to prove the following Krylov's type estimate for $X_t$: for any $f\in C^\infty\cap H^{-\alpha, p}$ and $T>0$, 
\begin{align}\label{Kr9}
\bE\left(\sup_{t\in[0,T]}\left|\int^t_0 f(X_s)\dif s\right|\right)\leq C\|f\|_{-\alpha,p}, 
\end{align}
where the constant $C$ is independent of $f$. 
In fact, if the above estimate is proven, then applying it to $b_n-b_m$, 
one sees that 
$(\int_0^\cdot b_n(X_s)\dif s)_{n\in\mN}$ is a Cauchy sequence in $L^1(\Omega; C([0,T]))$, and so there is a continuous 
adapted process denoted by $A^b_t$ such that
$$
\lim_{n\to\infty}\bE\left(\sup_{t\in[0,T]}\left|\int^t_0 b_n(X_s)\dif s-A^b_t\right|\right)=0.
$$
In order to show the above estimate, we need to have a better understanding for the following associated PDE
\begin{align}\label{EQ21}
\sL^a u-\lambda u+b\cdot\nabla u=f,
\end{align}
where $\lambda>0$,
$a^{ij}:=\sigma^{ik}\sigma^{jk}/2$ and $\sL^a u:=a^{ij}\p_i\p_j u$. 
Here and below, we use the usual Einstein's convention for summation: The same index appearing in 
a product will be summed automatically. 
In the sequel, in order to emphasize the dependence on $\sigma$, we sometimes write $\sL^\sigma:=\sL^a$.
Notice that the term $b\cdot\nabla u$ in \eqref{EQ21} should be understood in the distributional sense.

\medskip

Since the limiting process $t\mapsto A^b_t$ is usually not absolutely continuous (even not of finite variation) with respect to the Lebesgue measure, 
if there is no additional information of $A^b_t$,
it is in general hard to show the uniqueness, even in weak sense. 
In one dimensional case, 
when $b$ is the derivative of a $\gamma$-order H\"older continuous function with $\gamma\in(\frac{1}{2},1)$,
and $\sigma$ is $\frac{1}{2}$-order H\"older continuous and  bounded below by a positive constant, 
by using  the scaling function $s(x)=\int_0^x \exp \left(\int_0^y 2b(z)/\sigma^2(z)\dif z\right)\dif y$ to remove the drift
as well as Yamada-Watanabe's pathwise uniqueness result about one-dimensional SDE, 
Bass and Chen \cite{Ba-Ch1} showed the existence and uniqueness of strong solutions to SDE \eqref{SDE} in a special class of Dirichlet processes.
We also refer to \cite{Ch-En,Fl-Ru-Wo1,Fl-Ru-Wo2, Ru-Tr, He-Zh, Hu-Le-Mi} for more results about one dimensional SDEs driven by Brownian motion with distributional drifts.

\medskip 

However, in the multi-dimensional case, solving SDE \eqref{SDE} with singular  drift $b$ becomes quite involved. 
When $b\in L^p(\mR^d)$ for some $p>d$ and $\sigma=\mI_{d\times d}$, 
by Girsanov's transformation and $L^p$-theory of parabolic equations, 
Krylov and R\"ockner in \cite{Kr-Ro} showed that there is a unique strong solution to SDE \eqref{SDE}.
We also mention that the strong well-posedness of SDE \eqref{SDE} driven by multiplicative Brownian noise  
was studied in \cite{Xi-Zh, Zh0, Zh1}. 
Recently, when $\sigma\equiv\mI_{d\times d}$ and $b\in H^{-\alpha,p}$
with $\alpha\in(0,\frac{1}{2})$ and $p\in(\frac{d}{1-\alpha},\frac{d}{\alpha})$,
Flandoli, Issoglio and Russo \cite{Fl-Is-Ru} showed the existence and uniqueness of
``virtual''  solutions (a class of special weak solutions) to SDE \eqref{SDE}.
Let us make a brief introduction to their 
work. Denote
$\Phi(x):=x+\bu(x),$ 
where $\bu:\mR^d\to\mR^d$ 
solves PDE \eqref{EQ21} with $f=-b$ and $\sL^a=\frac{1}{2}\Delta$. 
For $\lambda$ being large enough, one can show that $\Phi$ is a $C^1$-diffeomorphism of $\mR^d$. 
Using It\^o's formula formally, it is easy to see that
$Y_t=\Phi(X_t)$ solves the following new SDE:
$$
Y_t=\Phi(x)+\int^t_0 \lambda \bu
\circ\Phi^{-1}(Y_s)\dif s+\int^t_0\nabla\Phi\circ\Phi^{-1}(Y_s)\dif B_s.
$$
Since this new SDE has continuous and non-degenerate diffusion coefficients, it is well known that the above SDE
admits a unique weak solution, and $X_t:=\Phi^{-1}(Y_t)$ is in turn defined as the solution of SDE \eqref{SDE} in \cite{Fl-Is-Ru}
(called ``virtual'' solution therein). The above $\Phi$ is usually called Zvonkin's transformation in literature (cf. \cite{Zv}).
It is noticed that the time-dependent drift $b$ is considered in \cite{Fl-Is-Ru} so that they need to solve a parabolic equation rather than 
an elliptic equation with distributional drift $b$. However, it is not answered whether
the above constructed $X$ really solves SDE \eqref{SDE} in the sense described in \eqref{SDE00}.
\medskip

The first purpose of this paper is to give an affirmative answer to the above question.
Now we outline the main points.
As mentioned above, in our work 
the crucial point is to prove the Krylov estimate \eqref{Kr9}.
Such an estimate together with the above transformation $\Phi$ will also lead to the weak uniqueness or the uniqueness of martingale solutions.
To achieve this aim, we need to tackle the following interesting problem: 
Find {\it minimal} conditions on $\Phi$ such that for some $C=C(\Phi,\alpha,p,d)>0$ and 
all distribution $f\in H^{-\alpha,p}$,
$$
\|f\circ\Phi\|_{-\alpha,p}\leq C\|f\|_{-\alpha,p}.
$$
Obviously this is a purely analytic problem, which has independent interest. In particular, the above estimate
implies that $T_\Phi(f):=f\circ\Phi$ is a bounded linear operator from $H^{-\alpha,p}$ to $H^{-\alpha,p}$.
We shall show it in Lemma \ref{Le21} below by using a duality argument. 
It should be  emphasized that our well-posedness result about SDE \eqref{SDE}
(see Theorem \ref{TH0} below) allows the drift $b$ being in the critical space $H^{-1/2,p}$.
Notice that this case is not covered in \cite{Ba-Ch1} and \cite{Fl-Is-Ru}, and which requires a
more delicate analysis for PDE \eqref{EQ21}.
Roughly to say, due to $b\in H^{-1/2,p}$, in order to make sense for $b\cdot\nabla u$, we need to at least assume $u\in H^{3/2,p}$.
Thus  $\sL^a u$ and $b\cdot\nabla u$ has the same order at scaling level. This is the source of the difficulty. 
Let us also mention that 
Bass and Chen  in \cite{Ba-Ch2} studied the weak well-posedness of SDE \eqref{SDE} when $b$ belongs to some generalized Kato's class, in particular, 
some measure-valued $b$ is allowed. 
Of course, our result is not comparable with \cite{Ba-Ch2}.

\medskip

The second aim of this paper is to show the existence and two-sided estimate of the heat kernel and 
the ergodicity associated with SDE \eqref{SDE}.
In fact, Zvonkin's transformation provides a satsifactory answer. 
In other words, if the transformed SDE admits a density and two-sided heat kernel  estimates, then the original SDE also admits a density and two-sided heat kernel  estimates.
Thus, one can construct the heat kernel of operator $\sL^a+b\cdot\nabla$ with distributional drift $b$. Notice
that when $b$ belongs to certain Kato's class,
the heat kernel of $\sL^a+b\cdot\nabla$ was constructed in \cite{ZhQ, Ch-Hu-Xi-Zh} by a perturbation argument. 
If it is not impossible, it seems hard to use the same perturbation method to study the heat kernel of $\sL^a+b\cdot\nabla$ when $b$ is a distribution.
Moreover, we also study the ergodicity of SDE \eqref{SDE} with $b=b^{(1)}+b^{(2)}$, where $b^{(1)}$ 
is the dissipative part and $b^{(2)}\in H^{-\alpha,p}$ is a distribution.
This is a continuation of work \cite{Xi-Zh}. Therein, when $b^{(2)}\in L^p$ for some $p>d$, the ergodicity is obtained by Zvonkin's transformation. It should be noticed that for the existence of invariant measures of SDE \eqref{SDE} with distributional drift $b$, 
a direct Lyapunov criterion (It\^o's formula) is not applicable.

\medskip

This paper is organized as follows: In Section 2, we prepare some analytic results.  In particular,  
product of two distributions in general Sobolev spaces is studied. 
In Section 3, we give the conceptions of martingale solutions and weak solutions, and prove their equivalence. 
In Section 4, we solve PDE \eqref{EQ21} with distributional drifts and variable coefficients by using Levi's freezing coefficient argument.
In Section 5, we state our main results and then prove them. 

\medskip

We close this section by mentioning 
some conventions used throughout this paper:
We use $:=$ as a way of definition. For $a, b\in \mR$, $a\vee b:= \max \{a, b\}$ and $a\wedge b:=\min \{a, b\}$,
and on $\mR^d$, $\nabla:=(\frac{\partial}{\partial x_1}, \dots, \frac{\partial}{\partial x_d})$ and 
$\Delta:= \sum_{k=1}^d \frac{\partial^2}{\partial x_k^2}$.
The letter $C$ with or without subscripts stands for an unimportant constant, whose value may change in difference places.
We use $A\asymp B$ to denote that $A$ and $B$ are comparable up to a constant, and use $A\lesssim B$ to denote $A\leq C B$ 
for some constant $C$.

\section{Preliminary}
In this section we present some analytic results that will be used later, and we believe that some of them have independent interest. 

Let $\rho$ be a nonnegative smooth function in $\mR^d$ with compact support in the unit ball and $\int\rho=1$. Define a family of mollifiers
$$
\rho_n(x)=n^d\rho(nx),\ \ n\in\mN.
$$
For a distribution $f\in\sD'$, if there is no further declaration, we always use $f_n$ to denote the mollifying approximation of $f$, that is,
$$
f_n(x):=f*\rho_n(x),
$$
where $*$ denotes the convolution in the distributional sense.
Let $\chi$ be a nonnegative smooth function with $\chi(x)=0$ for $|x|\geq 2$ and $\chi(x)=1$ for $|x|\leq 1$. For $R>0$, we shall also use the following cutoff function
\begin{align}\label{Cut}
\chi_R(x)=\chi(x/R).
\end{align}

\bd\label{def-Bessel}
For $\alpha\in\mR$ and $p\in[1,\infty)$, the Bessel potential space $H^{\a, p}$ is defined 
by
$$
H^{\alpha,p}:=(\mI-\Delta)^{-\alpha/2}(L^p)
$$
with norm
$$
\|f\|_{\alpha,p}:=\|(\mI-\Delta)^{\alpha/2} f\|_p,
$$
where $\|\cdot\|_p$ is the usual $L^p$-norm. We also denote by $H^{\alpha,p}_{loc}$ the space of all the distribution $f\in\sD'$ 
with $f\chi_R\in H^{\alpha,p}$ for any $R>0$, 
which is the local Bessel potential space.
\ed

For $\alpha\in(0,2)$ and $p\in(1,\infty)$, by Mihlin's multiplier theorem, we have
\begin{align}\label{Fr1}
\|f\|_{\alpha,p}\asymp \|(\mI-\Delta^{\alpha/2}) f\|_p\asymp \|f\|_p+\|\Delta^{\alpha/2}f\|_p,
\end{align}
where $\Delta^{\alpha/2}:=-(-\Delta)^{\alpha/2}$ is the usual fractional Laplacian, which has the following alternative expression up to a multiplying constant,
$$
\Delta^{\alpha/2}f(x)=\mbox{P.V.}\int_{\mR^d}\frac{f(x+y)-f(x)}{|y|^{d+\alpha}}\dif y,
$$
where P.V. stands for Cauchy's principle value. Clearly, if we write
$$
\Gamma_\alpha(f,g)(x):=\int_{\mR^d}(f(x+y)-f(x))(g(x+y)-g(x))\frac{\dif y}{|y|^{d+\alpha}},
$$
then
\begin{align}\label{Fr}
\Delta^{\alpha/2}(f g)=\Delta^{\alpha/2}f\cdot g+f\cdot\Delta^{\alpha/2}g+\Gamma_\alpha(f,g).
\end{align}
Notice that the following Sobolev's embedding holds:
\begin{align}\label{Fr2}
H^{\alpha,p}\subset 
\left\{
\begin{aligned}
&\cap_{q\in[p,dp/(d-p\alpha)]}L^q, &\mbox{ if } p\alpha<d,\\
&\cC^{\alpha-d/p}\cap(\cap_{q\geq p}L^q),\ &\mbox{ if } p\alpha>d,
\end{aligned}
\right. 
\end{align}
where $\cC^{\alpha-d/p}$ is the usual H\"older space.
Moreover, for any $\alpha\in(0,1]$ and $p\in(1,\infty)$, 
there is a constant $C=C(\alpha,p,d)>0$ such that for all $f\in H^{\alpha,p}$ (see \cite[Theorem 2.36]{Ba-Ch-Da}),
\begin{align}\label{GF1}
\|f(\cdot+y)-f(\cdot)\|_p\leq C |y|^\alpha\|\Delta^{\alpha/2} f\|_{p},
\end{align}
and if $p\alpha>d$, then for all $f\in H^{\alpha,p}$ and $x,y\in\mR^d$,
\begin{align}\label{RG1}
|f(x+y)-f(x)|\leq C |y|^{\alpha-\frac{d}{p}}\|\Delta^{\alpha/2} f\|_{p},
\end{align}
and the following Gagliardo-Nirenberg's inequality holds: for $p>1$ and $0<\alpha<\beta\leq 1$, and
all $f\in H^{\beta,p}\cap L^\infty$ (see \cite[Theorem 2.44]{Ba-Ch-Da}),
\begin{align}\label{Ga}
\|\Delta^{\alpha/2}f\|_{p\beta/\alpha} \leq C \|f\|_\infty^{1-\alpha/\beta} \|\Delta^{\beta/2} f \|_{p}^{\alpha/\beta}.
\end{align}

The following simple lemma plays a basic role in this paper.
\bl\label{Prod}
Let $p\in(1,\infty)$ and $\alpha\in(0,1]$ be fixed.
\begin{enumerate}[(i)]
\item For any $p_1,p_2\in[p,\infty)$ with $\frac{1}{p}\leq\frac{1}{p_1}+\frac{1}{p_2}<\frac{1}{p}+\frac{\alpha}{d}$, there is a constant $C>0$ such that
for all $f\in H^{\alpha,p_1}$ and $g\in H^{\alpha,p_2}$,
\be\label{EW1}
\|fg\|_{\alpha,p}\leq C \|f\|_{\alpha,p_1} \|g\|_{\alpha,p_2}.
\ee
In particular, if $p>d/\alpha$, then $H^{\alpha,p}$ is an algebra under pointwise product.
\item 
For any $p_1\in[p,\infty)$ 
and $p_2\in[\frac{p_1}{p_1-1},\infty)$ 
with $\frac{1}{p}\leq\frac{1}{p_1}+\frac{1}{p_2}<\frac{1}{p}+\frac{\alpha}{d}$, there is a constant $C>0$ such that
for all $f\in H^{-\alpha,p_1}$ and $g\in H^{\alpha,p_2}$,
\be\label{EW2}
\|fg\|_{-\alpha,p}\leq C \|f\|_{-\alpha,p_1} \|g\|_{\alpha,p_2}.
\ee
\end{enumerate}
\el
\bpf
(i) 
Below we fix $p_1,p_2\in[p,\infty)$ so that $\frac{1}{p}\leq\frac{1}{p_1}+\frac{1}{p_2}<\frac{1}{p}+\frac{\alpha}{d}$. 
By \eqref{Fr}, one sees that
$$
(\mI-\Delta^{\alpha/2})(fg)=(\mI-\Delta^{\alpha/2})f\cdot g-f\cdot\Delta^{\alpha/2}g-\Gamma_{\alpha}(f,g).
$$
Let $p'_i:=\frac{pp_i}{p_i-p}, i=1,2$ and 
$\delta:=\tfrac{d}{p'_1}-\tfrac{d}{p_2}+\alpha>0$.
By H\"older's inequality, \eqref{Fr2} and \eqref{GF1}, we have
\begin{align*}
\|fg\|_{\alpha,p}&\leq \|f\|_{\alpha,p_1}\|g\|_{p'_1}+\|f\|_{p'_2}\|g\|_{\alpha,p_2}+\int_{\mR^d}\|f(\cdot+y)-f(\cdot)\|_{p_1}\|g(\cdot+y)-g(\cdot)\|_{p'_1}\frac{\dif y}{|y|^{d+\alpha}}\\
&\lesssim \|f\|_{\alpha,p_1}\|g\|_{\alpha, p_2}+\|f\|_{\alpha,p_1}\|g\|_{\delta,p'_1}\int_{\mR^d}(1\wedge|y|^{\alpha+\delta})\frac{\dif y}{|y|^{d+\alpha}}
\lesssim \|f\|_{\alpha,p_1}\|g\|_{\alpha,p_2}.
\end{align*}
Thus we get \eqref{EW1}.
\\
\\
(ii) Let $p':=\frac{p}{p-1}$ and $p'_1:=\frac{p_1}{p_1-1}$. By the assumption, one sees that
$$
\tfrac{1}{p'_1}\leq \tfrac{1}{p'}+\tfrac{1}{p_2}<\tfrac{1}{p'_1}+\tfrac{\alpha}{d},\ p',p_2\in[p'_1,\infty).
$$
By duality and \eqref{EW1}, we have
\begin{align*}
\|fg\|_{-\alpha,p}&=\sup_{\|h\|_{p'}\leq 1}\left|\int_{\mR^d}fg\cdot (\mI-\Delta)^{-\alpha/2} h\right|
\leq \|f\|_{-\alpha, p_1}\sup_{\|h\|_{p'}\leq 1}\|g(\mI-\Delta)^{-\alpha/2} h\|_{\alpha,p_1'}\\
&\lesssim\|f\|_{-\alpha, p_1} \|g\|_{\alpha,p_2} \sup_{\|h\|_{p'}\leq 1} \|(\mI-\Delta)^{-\alpha/2} h\|_{\alpha,p'}
= \|f\|_{-\alpha, p_1}\|g\|_{\alpha,p_2}.
\end{align*}
Thus we get \eqref{EW2}.
\epf
\br\label{Re22}\rm
(i) By the above lemma, one sees that if $f\in H^{\pm\alpha,p_1}_{loc}$ and $g\in H^{\alpha,p_2}_{loc}$ 
with $p, p_1, p_2$ being as in the lemma,
then $fg\in H^{\pm\alpha,p}_{loc}$.
Moreover, from the proof of the lemma and by \eqref{Ga}, we also have
\begin{align}
\|fg\|_{-\alpha,p}\leq C\|f\|_{-\alpha,p}(\|g\|_\infty+\|\Delta^{\alpha/2}g\|_{p_2}),\ \ p_2>\tfrac{p}{p-1}\vee\tfrac{d}{\alpha}.
\end{align}
(ii) Let $\alpha\in\mR$.
For H\"older-Besov space $\cC^\alpha$, it is well known that for $\beta>\alpha>0$ (for example, see \cite{Ba-Ch-Da} and \cite{Gu-Im-Pe}),
$$
\|fg\|_{\cC^{-\alpha}}\leq C\|f\|_{\cC^{-\alpha}}\|g\|_{\cC^\beta}.
$$
Compared with \eqref{EW2}, the duality argument can not be used to show the above inequality for $\beta=\alpha$ since 
the dual space of $\cC^{-\alpha}$ does not equal $\cC^{\alpha}$. 
\er
Let $\cD^0_\infty$ be the set of all  $C^1$-diffeomorphisms on $\mR^d$:
$$
\cD^0_\infty:=\Big\{\Phi:\mR^d\to\mR^d,\ \ \|\Phi\|_{\cD^0_\infty}:=\|\nabla\Phi\|_\infty+\|\nabla\Phi^{-1}\|_\infty<\infty\Big\}.
$$ 
Clearly, $\cD^0_\infty$ is closed under the inverse operation, that is, $\Phi\in \cD^0_\infty$ implies $\Phi^{-1}\in \cD^0_\infty$.
The following lemma is easy by interpolation.

\bl\label{Le211}Let $\Phi\in\cD^0_\infty$ be a $C^1$-diffeomorphism. For any $\alpha\in[0,1]$ and $p>1$, there is a constant
$C=C(\alpha,d,p,\|\Phi\|_{\cD^0_\infty})>0$ such that for all $f\in H^{\alpha,p}$,
\begin{align}\label{KJ1}
\|f\circ\Phi\|_{\alpha,p}\leq C\|f\|_{\alpha,p}.
\end{align}
\el
\begin{proof}
By the change of variable, it is clear that
\begin{align}\label{Fr5}
\|f\circ\Phi\|_p\leq \|\det(\nabla\Phi^{-1})\|_\infty^{1/p}\|f\|_p.
\end{align}
On the other hand, noticing that 
\begin{align*}
\nabla (f\circ\Phi)=(\nabla f)\circ\Phi\cdot\nabla\Phi,
\end{align*}
we have
$$
\|\nabla (f\circ\Phi)\|_p\leq\|(\nabla f)\circ\Phi\|_p\cdot\|\nabla\Phi\|_\infty\leq  \|\det(\nabla\Phi^{-1})\|_\infty^{1/p}\|\nabla\Phi\|_\infty\|\nabla f\|_p.
$$
Hence, \eqref{KJ1} holds for $\alpha=1$. 
Noticing that
$T_\Phi: f\mapsto f\circ \Phi$ is linear, by the interpolation theorem, we get \eqref{KJ1} for $\alpha\in[0,1]$.
\end{proof}

To perform Zvonkin's transformation, we need to show that \eqref{KJ1} holds for {\it negative} $\alpha$.
To this aim, we introduce a subclass of $\cD^0_\infty$ as follows: For $\beta\in(0,1]$ and $q\in(d/\beta,\infty)$, 
$$
\cD^\beta_q:=\Big\{\Phi\in\cD^0_\infty: \|\Phi\|_{\cD^\beta_q}:=\|\Phi\|_{\cD^0_\infty}+\|\mI-\nabla\Phi\|_{\beta,q}<\infty\Big\}.
$$
The following proposition shows that $\cD^\beta_q$ is closed under the inverse operation.
\bp\label{Pr22}
Let $\beta\in(0,1]$ and $q\in(d/\beta,\infty)$. For any $\Phi\in\cD^\beta_q$, we have $\Phi^{-1}\in\cD^\beta_q$ and
$$
\|\det(\nabla\Phi)-1\|_{\beta,q}, \|\det(\nabla\Phi^{-1})-1\|_{\beta,q}<\infty.
$$
\ep
\begin{proof}
(i) Let $U(x):=\nabla\Phi(x)-\mI$. By the definition of the determinant of a matrix, one sees that
\begin{align}\label{NL1}
\det(\nabla\Phi)=\det(\mI+U)=1+P(U),
\end{align}
where $P$ is a polynomial of $(U_{ij})$ without zero order term. Due to $q>d/\beta$, by \eqref{EW1} with $p_1=p_2=q$, 
we have
\begin{align}\label{NL2}
U\in H^{\beta,q}\Rightarrow U^m\in H^{\beta,q}\mbox{ for any $m\in\mN$}\Rightarrow P(U)\in H^{\beta,q}.
\end{align}
Hence,
$$
\|\det(\nabla\Phi)-1\|_{\beta,q}<\infty.
$$
(ii) To prove $\Phi^{-1}\in\cD^\beta_q$, by definition it suffices to show $\|\mI-\nabla\Phi^{-1}\|_{\beta,q}<\infty$. 
First of all, since $\nabla\Phi^{-1}=(\nabla\Phi)^{-1}\circ\Phi^{-1}$, by \eqref{KJ1} we have
$$
\|\mI-\nabla\Phi^{-1}\|_{\beta,q}\lesssim\|\mI-(\nabla\Phi)^{-1}\|_{\beta,q}
=\|(\nabla\Phi)^{-1}(\nabla\Phi-\mI)\|_{\beta,q}=\|(\nabla\Phi)^{-1}U\|_{\beta,q}.
$$
Clearly,
$$
\|(\nabla\Phi)^{-1}U\|_{q}\leq\|(\nabla\Phi)^{-1}\|_\infty\|U\|_q<\infty. 
$$
Noticing that
$$
(\nabla\Phi)^{-1}(x+z)-(\nabla\Phi)^{-1}(x)=(\nabla\Phi)^{-1} (x)\left( U(x)-U(x+z) \right)(\nabla\Phi)^{-1}(x+z),
$$
by definition \eqref{Fr}, we can write
\begin{align*}
&\Delta^{\beta/2}((\nabla\Phi)^{-1}U)(x)=\int_{\mR^d} 
\left[ (\nabla\Phi)^{-1}(x+z)-(\nabla\Phi)^{-1}(x)\right]U(x+z) \frac{\dif z}{|z|^{d+\beta}}+\Big((\nabla\Phi)^{-1}\Delta^{\beta/2}U\Big)(x)\\
&\qquad=\int_{\mR^d} 
(\nabla\Phi)^{-1} (x)\left( U(x)-U(x+z) \right)((\nabla\Phi)^{-1}U)(x+z) \frac{\dif z}{|z|^{d+\beta}}+\Big((\nabla\Phi)^{-1}\Delta^{\beta/2}U\Big)(x). 
\end{align*}
Let $K(x,z):=((\nabla\Phi)^{-1}U)(x+z)$. 
Since by \eqref{Fr2} and \eqref{NL2}, 
$$
\| \nabla \Phi\|_{\cC^{\beta-d/q}}\lesssim \| \nabla \Phi\|_{\beta,q}  \leq \|U\|_{\beta, q}+1<\infty,
$$
it is easy to see that
$$
\|K(\cdot,0)\|_\infty<\infty,\ |K(x,z)-K(x,0)|\lesssim (1\wedge|z|)^{\beta-d/q}.
$$
Therefore,  
by \eqref{GF1},
\begin{align*}
\left\|\Delta^{\beta/2}((\nabla\Phi)^{-1}U)\right\|_q&\leq
\left\|\int_{\mR^d} (\nabla\Phi)^{-1}(\cdot)\left( U(\cdot)-U(\cdot+z) \right)(K(\cdot,z)-K(\cdot,0))\frac{\dif z}{|z|^{d+\beta}}\right\|_q\\
&\quad+\left\|(\nabla\Phi)^{-1}(\Delta^{\beta/2}U)K(\cdot,0)\right\|_q+\left\|(\nabla\Phi)^{-1}\Delta^{\beta/2}U\right\|_q\\
&\lesssim\|(\nabla\Phi)^{-1}\|_\infty\int_{\mR^d}\|U(\cdot)-U(\cdot+z)\|_q(1\wedge|z|)^{\beta-d/q}\frac{\dif z}{|z|^{d+\beta}}\\
&\quad+\|(\nabla\Phi)^{-1}\|_\infty\|\Delta^{\beta/2}U\|_q\|K(\cdot,0)\|_\infty+\|(\nabla\Phi)^{-1}\|_\infty\|\Delta^{\beta/2}U\|_q\\
&\lesssim \|U\|_{\beta,q}\int_{\mR^d}(1\wedge|z|)^{2\beta-d/q}\frac{\dif z}{|z|^{d+\beta}}+\|U\|_{\beta,q}<\infty.
\end{align*}
The proof is complete. 
\end{proof}

Now we can show an analogous version of \eqref{KJ1} for $\a<0$,  
which is crucial for applying Zovnkin's transformation. 
\bl\label{Le21}
Let $\Phi\in\cD^\beta_q$ for some $\beta\in(0,1)$ and $q\in(d/\beta,\infty)$.  For any $\alpha\in[0,\beta]$ and $p>\frac{d}{d-\alpha}$, there is a constant 
$C=C(\alpha,\beta,d,p,\|\Phi\|_{\cD^\beta_q})>0$ 
such that for all $f\in H^{-\alpha,p}$,
\begin{align}\label{KJ2}
\|f\circ\Phi\|_{-\alpha,p}\leq C\|f\|_{-\alpha,p}.
\end{align}
\el
\begin{proof}
By a density argument, we may assume $f\in \sD$. Letting $p':=p/(p-1)$, we have
\begin{align*}
\|f\circ\Phi\|_{-\alpha,p}&=\sup_{\|g\|_{p'}\leq 1}\left|\int_{\mR^d}f\circ\Phi(x)\cdot (\mI-\Delta)^{-\alpha/2}g(x)\dif x\right|\\
&=\sup_{\|g\|_{p'}\leq 1}\left|\int_{\mR^d}f(x)\cdot (\mI-\Delta)^{-\alpha/2}g\circ\Phi^{-1}(x)\cdot \left|\det(\nabla\Phi^{-1}(x))\right| \dif x\right|\\
&\leq \|f\|_{-\alpha,p}\sup_{\|g\|_{p'}\leq 1}\left\|(\mI-\Delta)^{-\alpha/2}g\circ\Phi^{-1}\cdot\det(\nabla\Phi^{-1})\right\|_{\alpha,p'}.
\end{align*}
Write $G_\alpha:=(\mI-\Delta)^{-\alpha/2}g\circ\Phi^{-1}$. 
Since $p>\frac{d}{d-\alpha}$ and $q>\frac{d}{\beta}$, we can choose $p_2\geq p'$ so that
$$
\tfrac{1}{q}-\tfrac{\beta-\alpha}{d}\leq\tfrac{1}{p_2}<\tfrac{\alpha}{d}<\tfrac{1}{p'}.
$$
Thus, by \eqref{EW1} with $p_1=p'$ and the above $p_2$, 
\begin{align*}
\|f\circ\Phi\|_{-\alpha,p}&\lesssim \|f\|_{-\alpha,p}\sup_{\|g\|_{p'}\leq 1}\Big(\left\|G_\alpha\cdot(\det(\nabla\Phi^{-1})-1)\right\|_{\alpha,p'}
+\|G_\alpha\|_{\alpha,p'}\Big)\\
&\lesssim \|f\|_{-\alpha,p}\sup_{\|g\|_{p'}\leq 1}\|G_\alpha\|_{\alpha,p'}\Big(\|\det(\nabla\Phi^{-1})-1)\|_{\alpha,p_2}+1\Big)\\
&\!\!\!\!\!\!\!\stackrel{\eqref{KJ1},
\eqref{Fr2}}{\lesssim} \|f\|_{-\alpha,p}\Big(\|\det(\nabla\Phi^{-1})-1)\|_{\beta,q}+1\Big),
\end{align*}
which gives \eqref{KJ2} since $\|\det(\nabla\Phi^{-1})-1)\|_{\beta,q}$ is finite by Proposition \ref{Pr22}.
\end{proof}
\br\label{Re35}\rm
By estimate \eqref{KJ2},
 for any $\Phi\in \cD^{\beta}_q$ and  $f\in H^{-\alpha, p}_{loc}$ with $\beta, q$ and $\alpha,p$ being as in the above lemma, 
 we can define a distribution $f\circ \Phi$ by 
$$
\<f\circ \Phi,g\>= \<f,  g\circ\Phi^{-1}\cdot |\det( \nabla \Phi^{-1})|\>,\ \ g\in \sD.
$$ 
In particular, it makes sense that $(f\circ\Phi^{-1})\circ\Phi=f$  for any $f\in H^{-\alpha,p}_{loc}$. 
\er

\section{Martingale problems and weak solutions}
Let $\mC$ be the space of all continuous functions from $\mR_+$ to $\mR^d$, which is endowed with the usual Borel $\sigma$-field $\cB(\mC)$.
All the probability measures over $(\mC,\cB(\mC))$ is denoted by $\sP(\mC)$. 
Let $w_t$ be the coordinate process over $\mC$, that is,
$$
w_t(\omega)=\omega_t,\ \ \omega\in\mC.
$$
For $t\geq 0$, let $\cB_t(\mC)$ be the natural filtration generated by $\{w_s: s\leq t\}$.
For given $R>0$, we shall use the following truncated $\cB_t(\mC)$-stopping time
\begin{align}\label{Tau}
\tau_R:=\inf\{t>0: |w_t|>R\}.
\end{align}
Notice that for each $\omega\in\mC$, it automatically holds that
\begin{align}\label{DF1}
\lim_{R\to\infty}\tau_R(\omega)=\infty.
\end{align}
For a probability measure $\mP\in\sP(\mC)$, the expectation with respect to $\mP$ will be denoted by $\mE^\mP$ or simply by $\mE$
if there is no confusion. 

Now we introduce the following important notion for later use.
\bd
(Local Krylov's estimate) Let $\alpha\in[0,1]$ and $p>1$.
We call a probability measure $\mP\in\sP(\mC)$ satisfy local Krylov's estimate with indices $\alpha,p$ if for any $T>0$ and $R\geq 1$, 
there are positive constants $C_{T,R}$ and $\gamma$ 
such that for all $f\in C^\infty$, $0\leq t_0<t_1\leq T$ and $\tau\leq\tau_R$,
\begin{align}\label{Kr}
\mE\left|\int^{t_1\wedge\tau}_{t_0\wedge\tau}f(w_s)\dif s\right|^2\leq C_{T,R}(t_1-t_0)^{1+\gamma}\|f\chi_R\|^2_{-\alpha,p}.
\end{align}
If $C_{T,R}$ does not depend on $R$, then the above estimate will be called global Krylov's estimate.
All the probability measure $\mP$ with property \eqref{Kr} is denoted by $\sK^{\alpha}_{p}(\mC)$.
\ed

About this definition we have the following useful consequence.
\bp\label{Pr13}
Let $\alpha\in[0,1]$, $p>1$  and $\mP\in\sK^\alpha_p(\mC)$. For any $f\in H^{-\alpha,p}_{loc}$, 
there is a continuous $\cB_t(\mC)$-adapted process $A^f_t$ such that 
for any mollifying approximation $f_n=f*\rho_n$ and any $T>0$, 
\begin{align}\label{NL3}
\lim_{n\to\infty}\mE\left(\sup_{t\in[0,T]}\left|\int^t_0f_n(w_s)\dif s-A^f_t\right|\wedge 1\right)=0.
\end{align}
Moreover, for each $R\geq 1$, the mapping $H^{-\alpha,p}\ni f\mapsto A^f_{\cdot\wedge\tau_R}\in L^2(\mC,\mP; C([0,T]))$ is 
a bounded linear operator, where $\tau_R$ is defined in \eqref{Tau}, and for all $0\leq t_0<t_1\leq T$,
\begin{align}\label{Kr0}
\mE\left|A^f_{t_1\wedge\tau_R}-A^f_{t_0\wedge\tau_R}\right|^2\leq C_{T,R}(t_1-t_0)^{1+\gamma}\|f\chi_R\|^2_{-\alpha,p},
\end{align}
where the constants $C_{T,R}$ and $\gamma$ are the same as in \eqref{Kr}.
\ep
\bpf
Let $R\geq 1$. For any $f\in C^\infty$ and $T>0$, by \eqref{Kr} and Kolmogorov's continuity criterion (see \cite{Re-Yo}), we have
$$
\mE\left(\sup_{t\in[0,T]}\left|\int^{t\wedge\tau_R}_0f(w_s)\dif s\right|^2\right)\leq C_{T,R}\|f\chi_R\|_{-\alpha,p}^2.
$$
In particular, applying this to smooth function $f_n-f_m$, we get
\begin{align*}
&\lim_{n,m\to\infty}\mE\left(\sup_{t\in[0,T]}\left|\int^{t\wedge\tau_R}_0(f_n-f_m)(w_s)\dif s\right|^2\right)
\leq C_{T,R}\lim_{n,m\to\infty}\|(f_n-f_m)\chi_R\|^2_{-\alpha,p}\\
&\quad=C_{T,R}\lim_{n,m\to\infty}\|((f\chi_{2R})_n-(f\chi_{2R})_m)\chi_R\|^2_{-\alpha,p}
\stackrel{\eqref{EW2}}{\leq} C'_{T,R}\lim_{n,m\to\infty}\|(f\chi_{2R})_n-(f\chi_{2R})_m\|^2_{-\alpha,p}=0,
\end{align*}
which means that 
$(\int^{\cdot\wedge\tau_R}_0f_n(w_s)\dif s)_{n\in\mN}$ 
is a Cauchy sequence in $L^2(\mC,\mP; C([0,T]))$. So, for each $R\geq 1$,  
there is a continuous $\cB_t(\mC)$-adapted process  $A^{f,R}_t$ such that for all $T>0$,
$$
\lim_{n\to\infty}\mE\left(\sup_{t\in[0,T]}\left|\int^{t\wedge\tau_R}_0f_n(w_s)\dif s-A^{f,R}_t\right|^2\right)=0.
$$
In particular, there is a $\mP$-null set $N$ such that for all $\omega\notin N$ and 
$R,R'\in\mN$ with $R<R'$,
$$
A^{f,R}_t(\omega)=A^{f,R'}_t(\omega),\ \forall t<\tau_R(\omega).
$$
Since $\lim_{R\to\infty}\tau_R(\omega)=\infty$ (see \eqref{DF1}), 
we may define a continuos adapted process $A^f_t$ on $\mR_+$ by
$$
A^f_t(\omega):=A^{f,R}_t(\omega),\ \ t\leq\tau_R(\omega),\ \omega\notin N.
$$
Now for any $\eps>0$, by Chebyshev's inequality we have
\begin{align*}
\mP\left(\sup_{t\in[0,T]}\left|\int^t_0f_n(w_s)\dif s-A^f_t\right|>\eps\right)&\leq
\mP\left(\sup_{t\in[0,T]}\left|\int^{t\wedge\tau_R}_0f_n(w_s)\dif s-A^f_t\right|>\eps\right)+\mP(\tau_R\leq T)\\
&\leq\mE\left(\sup_{t\in[0,T]}\left|\int^{t\wedge\tau_R}_0f_n(w_s)\dif s-A^f_t\right|^2\right)/\eps^2+\mP(\tau_R\leq T),
\end{align*}
which converges to zero by first letting $n\to\infty$ and then $R\to\infty$ and \eqref{DF1}. Thus, we get \eqref{NL3}. As for \eqref{Kr0}, it follows by \eqref{Kr}.
\epf
\br\label{Re23}
\rm
\begin{enumerate}[(i)]
\item Estimate \eqref{Kr0} implies that $t\mapsto A^f_t$ is a locally zero energy process, that is, for any $R\geq 1$,
$$
\lim_{\delta\to 0}\sup_{\{\Pi_t: \rm{mesh}(\Pi_t)<\delta\}}\sum_{i=0}^n\mE|A^f_{t_{i+1}\wedge\tau_R}-A^f_{t_{i}\wedge\tau_R}|^2=0,
$$
where $\Pi_t:=\{t_0,t_1,\cdots,t_n\}$ denotes a partition of $[0,t]$.
\item
If $f\in L^q_{loc}(\mR^d)$ with $q\geq pd/(d+p\alpha)$, 
then $t\mapsto A^f_t$ is absolutely continuous and
$$
A^f_t=\int^t_0f(w_s)\dif s.
$$
Indeed, it follows by Sobolev's embedding $L^q_{loc}\subset H^{-\alpha,p}_{loc}$.
\end{enumerate}
\er
\bd\label{Def2}
(Martingale Problem) Let $\alpha\in[0,1]$ and $p>1$.
We call a probability measure $\mP\in\sK^\alpha_p(\mC)$ a martingale solution of SDE \eqref{SDE} with starting point $x\in\mR^d$
if for any $f\in C^\infty$,
\begin{align}\label{EQ1}
M^f_t:=f(w_{t})-f(x)-\int^{t}_0(\sL^\sigma f)(w_s)\dif s-A^{b\cdot\nabla f}_{t}
\end{align}
is a continuous local $\cB_t(\mC)$-martingale with $M^f_0=0$ under $\mP$, provided that $b\cdot \nabla f\in H^{-\a,p}_{loc}$, 
where $\sL^\sigma f:=\sigma^{ik}\sigma^{jk}\p_i\p_j f/2$.
All the martingale solution $\mP\in\sK^\alpha_p(\mC)$ of SDE \eqref{SDE} with coefficients $\sigma,b$ and starting point $x$ 
is denoted by $\sM^{\alpha,p}_{\sigma,b}(x)$.
\ed
\br\rm
\begin{enumerate}[(i)]
\item In the above definition, in order to make $A^{b\cdot\nabla f}_t$ well defined, 
we need to at least assume $b\in H^{-\alpha,p}_{loc}$ by Proposition \ref{Pr13} and Lemma \ref{Prod}.
The localization sequence of stopping times for $M^f_t$ being a martingale can be taken as $\tau_n$ (see \eqref{Tau}).
Moreover, for $\mP\in \sM^{\a,p}_{\sigma,b}(x)$,  due to $M^f_0=0$, we have $\mP(w_0=x)=1$. 

\item Trivially $\sM^{0,\infty}_{\sigma,b}(x)$ is the usual notion of martingale solutions (see \cite{St-Va}).
\end{enumerate}
\er

As a direct consequence of martingale solutions, we have
\bl[Generalized It\^o's formula]\label{Re1}Let $\alpha\in(0,\frac{1}{2}]$, $p>\frac{d}{1-\alpha}$ and $\beta\in[\alpha,1]$, $q\in(\frac{d}{\beta},\infty)$. 
Suppose $\sigma\in H^{\beta,q}_{loc}$ and $b\in H^{-\alpha,p}_{loc}$. 
For any $f\in H^{2-\alpha,p}_{loc}$ and $\mP\in\sM^{\alpha,p}_{\sigma,b}(x)$,  
$$
M^f_t:=f(w_t)-f(x)-
A^{(\sL^\sigma+b\cdot\nabla) f}_t
$$ 
is a continuous local $\cB_t(\mC)$-martingale under $\mP$.
\el
\bpf
Let $f_n:=f*\rho_n$ be the mollifying approximation of $f$. Fix $R>0$. By Definition \ref{Def2},
the process $t\mapsto M^{f_n}_{t\wedge\tau_R}$ is a continuous martingale.
Since $\mP(\lim_{R\to\infty}\tau_R=\infty)=1$, 
it suffices to show
$$
\lim_{n\to\infty}\mE\left(\sup_{t\in[0,T]}|M^{f_n}_{t\wedge\tau_R}-M^{f}_{t\wedge\tau_R}|^2\right)
\leq 4\lim_{n\to\infty}\mE|M^{f_n}_{T\wedge\tau_R}-M^{f}_{T\wedge\tau_R}|^2=0,\ \ R>0,
$$
where the first inequality is due to Doob's maximal inequality. 
By \eqref{Kr0}, we only need to show
$$
\lim_{n\to\infty}\|(\sL^\sigma+b\cdot\nabla) (f_n-f)\cdot\chi_R\|_{-\alpha,p}=0.
$$
Since $f\in H^{2-\alpha,p}_{loc}$, for any $R\geq 1$, we have $\nabla^2 (f \chi_R) \in H^{-\alpha,p}$.   
Noticing $\sigma\in H^{\beta,q}_{loc}$, 
by \eqref{EW2} with $p_1=p$ and $p_2=qd/(d-q(\beta-\alpha))>d/\alpha>p/(p-1)$, one sees that 
\begin{align*}
&\|(\sL^\sigma (f_n-f))\chi_R\|_{-\alpha,p}=\|(\sL^{\sigma\chi_{2R}} ((f_n-f)\chi_{2R}))\chi_R\|_{-\alpha,p}\\
&\lesssim\|\sL^{\sigma\chi_{2R}} ((f_n-f)\chi_{2R})\|_{-\alpha,p}
\lesssim\|\nabla^2((f_n-f)\chi_{2R})\|_{-\alpha,p}\|(\sigma\chi_{2R})(\sigma\chi_{2R})^*\|_{\alpha,p_2}\\
&\stackrel{\eqref{Fr2}}{\lesssim} \|\nabla^2((f_n-f)\chi_{2R})\|_{-\alpha,p}\lesssim \|(f_n-f)\chi_{2R}\|_{2-\alpha,p}\to 0 \ \mbox{ as $n\to\infty$.}
\end{align*}
Similarly, for $p_1=p$ and $p_2=pd/(d-p(1-2\alpha))>d/\alpha$, we have
\begin{align*}
\|b\cdot\nabla (f_n-f)\cdot\chi_R\|_{-\alpha,p}
&=\|b\chi_R\cdot\nabla ((f_n-f)\chi_{2R})\|_{-\alpha,p}\stackrel{\eqref{EW2}}{\lesssim} \|b\chi_R\|_{-\alpha,p}\|\nabla ((f_n-f)\chi_{2R})\|_{\alpha,p_2}\\
&\lesssim\|b\chi_R\|_{-\alpha,p}\|(f_n-f)\chi_{2R}\|_{1+\alpha,p_2}\stackrel{\eqref{Fr2}}{\lesssim}
\|b\chi_R\|_{-\alpha,p}\|(f_n-f)\chi_{2R}\|_{2-\alpha,p}\to 0.
\end{align*}
The proof is complete.
\epf

\bp\label{Pr24}
(Zvonkin's transformation) 
Let $\alpha\in(0,\frac{1}{2}]$, $p>\frac{d}{1-\alpha}$ and $\beta\in[\alpha,1]$, $q\in(\frac{d}{\beta},\infty)$.
Suppose that $\sigma\in H^{\beta,q}_{loc}$,  $b\in H^{-\alpha,p}_{loc}$ and $\Phi\in\cD^{1-\alpha}_p$.
Define 
\be\label{def-new}
\tilde\sigma:=(\nabla\Phi\cdot\sigma)\circ\Phi^{-1}, \ \ \tilde b:=(\sL^\sigma\Phi+b\cdot\nabla\Phi)\circ\Phi^{-1}.
\ee
Then we have

\begin{enumerate}[(i)] 
\item $\tilde b\in H^{-\alpha,p}_{loc}$ and $\tilde\sigma\in H^{\beta',q'}_{loc}$ for $\beta':=\beta\wedge(1-\alpha)$ and
\begin{align}\label{GS1}
\tfrac{1}{q'}:=\left\{
\begin{aligned}
&\tfrac{1}{q}\vee \left(\tfrac{1}{p}- \tfrac{1-\a-\beta}{d}\right),\ \ \beta\in[\alpha,1-\alpha],\\
&\tfrac{1}{p}\vee \left(\tfrac{1}{q}-\tfrac{\alpha+\beta-1}{d}\right),\ \ \beta\in(1-\alpha,1], 
\end{aligned}
\right. 
\end{align}
and also $q'>d/\beta'$. 

\item For any $x\in\mR^d$, it holds that
\begin{align}\label{Zvon}
\mP\in\sM^{\alpha,p}_{\sigma,b}(x)\Leftrightarrow
\mP\circ\Phi^{-1}\in\sM^{\alpha,p}_{\tilde\sigma,\tilde b}(\Phi(x)).
\end{align}
Here $\mP\circ\Phi^{-1}$ means that for 
$A\in\cB(\mC)$, $\mP\circ\Phi^{-1}(A)=\mP\left(\{\omega: \Phi(w_\cdot(\omega))\in A\}\right)$.
\end{enumerate}
\ep
\begin{proof}
(i) It follows by Remark \ref{Re22} (i), Lemmas \ref{Le211} and \ref{Le21}.
\\
\\
(ii)  Since $\beta'\in[\alpha,1)$ and $q'>d/\beta'$, by symmetry, we only show $\Rightarrow$.
To show that $\mP\circ\Phi^{-1}$ is a martingale solution of SDE \eqref{SDE} with coefficients $\tilde\sigma$ and $\tilde b$, one only needs to check that
$\mP\circ\Phi^{-1}\in\sK^\alpha_p(\mC)$ and  for any $f\in C^\infty$,
$$
\tilde M^f_t:=f(w_t)-f(w_0)-\int^t_0\sL^{\tilde \sigma}f(w_s)\dif s-A^{\tilde b\cdot\nabla f}_t
$$
is a continuous local $\cB_t(\mC)$-martingale under $\mP\circ\Phi^{-1}$.
First of all, since $\Phi$ is a homeomorphism, there is an $R'>R$ large enough so that for any $\tau\leq\tau_R$,
$$
\tau\circ\Phi\leq\tau_R\circ\Phi\leq \tau_{R'}, 
$$
and $\tau\circ \Phi$ is also a $\cB_t(\mC)$-stopping time. 
Thus,  for any $\tau\leq\tau_R$ and
$0\leq t_0<t_1\leq T$,
we have
\begin{align*}
&\mE^{\mP\circ\Phi^{-1}}\left|\int^{t_1\wedge\tau}_{t_0\wedge\tau}f(w_s)\dif s\right|^2
=\mE\left|\int^{t_1\wedge\tau\circ\Phi}_{t_0\wedge\tau\circ\Phi}(f\chi_R)\circ\Phi (w_s)\dif s\right|^2\\
&\quad\stackrel{\eqref{Kr}}{\leq} 
C_{R',T}(t_1-t_0)^{1+\gamma}\|(f\chi_R)\circ\Phi \cdot\chi_{R'}\|^2_{-\alpha,p}\stackrel{\eqref{EW1},\eqref{KJ2}}{\leq} 
C_{R',T}(t_1-t_0)^{1+\gamma}\|f\chi_R\|^2_{-\alpha, p}.
\end{align*}

Next we show that $\tilde M^f_t$
is a continuous local $\cB_t(\mC)$-martingale under $\mP\circ\Phi^{-1}$.
By definition, it suffices to prove that $\tilde M^f_t\circ\Phi$ is a continuous local $\cB_t(\mC)$-martingale under $\mP$. 
Noticing that
$$
\nabla (f\circ\Phi)=(\nabla f)\circ\Phi\cdot\nabla\Phi,
$$
and  in the distributional sense, 
$$
\p^2_{ij}(f\circ\Phi)=(\p_k f)\circ\Phi\cdot\p^2_{ij}\Phi^k+(\p^2_{i'j'} f)\circ\Phi\cdot\p_i\Phi^{i'}\cdot\p_j\Phi^{j'},
$$
we have
\begin{align}\label{YP1}
(\sL^{\tilde \sigma}f)\circ\Phi=\sL^\sigma(f\circ\Phi)-\sL^\sigma\Phi\cdot\nabla f\circ\Phi,
\end{align}
and by Remark \ref{Re35}, 
\begin{align}\label{YP2}
(\tilde b\cdot\nabla f)\circ\Phi=\sL^\sigma\Phi\cdot\nabla f\circ\Phi+b\cdot\nabla(f\circ\Phi).
\end{align}
Hence,
\begin{align*}
\tilde M^f_t\circ\Phi&=f\circ\Phi(w_t)-f\circ\Phi(w_0)-\int^t_0(\sL^{\tilde \sigma}f)\circ\Phi(w_s)\dif s-A^{(\tilde b\cdot\nabla f)\circ\Phi}_t\\
&=f\circ\Phi(w_t)-f\circ\Phi(w_0)-A^{\sL^\sigma(f\circ\Phi)}_t-A^{b\cdot\nabla(f\circ\Phi)}_t. 
\end{align*}
Moreover, since by \eqref{EW1} and \eqref{KJ1},
\begin{align*}
&\|\nabla (f\circ\Phi)\cdot\chi_R\|_{1-\alpha,p}=\|(\nabla f)\circ\Phi\cdot\nabla\Phi\cdot\chi_R\|_{1-\alpha,p}\\
&\leq\|((\nabla f)\circ\Phi)\cdot\chi_R\cdot(\nabla\Phi-\mI)\|_{1-\alpha,p}
+\|(\nabla f)\circ\Phi\cdot\chi_R\|_{1-\alpha,p}\\
&\lesssim\|(\nabla f)\circ\Phi\cdot\chi_R\|_{1-\alpha,p}\cdot\|\nabla\Phi-\mI\|_{1-\alpha,p}+\|\nabla f\cdot(\chi_R\circ\Phi^{-1})\|_{1-\alpha,p}<\infty,
\end{align*}
we have $f\circ\Phi\in H^{2-\alpha,p}_{loc}$. By Lemma \ref{Re1},
$\tilde M^f_t\circ\Phi=M^{f\circ\Phi}_t$ is a continuous local martingale with respect to $\mP$. 
\end{proof}

\br\rm
The importance of  \eqref{Zvon} lies in the fact that if there is one and only one element in $\sM^{\alpha,p}_{\tilde\sigma,\tilde b}(\Phi(x))$,
then there is automatically one and only one element  in $\sM^{\alpha,p}_{\sigma,b}(x)$.
Moreover, the heat kernel estimates and ergodicity can also be derived by \eqref{Zvon}. 
\er

 
Next we introduce the notion of weak solutions and discuss the relationship between martingale solutions and weak solutions.

\bd[Weak solutions]\label{Def99}
 Let $\sigma$ be locally bounded and $b\in H^{-\a, p}_{loc}$ for some $\a\in [0,1]$ and $p>1$.
 Let $(X,B)$ be two $\mR^d$-valued continuous adapted processes on some 
 filtered probability space $(\Omega,\sF, (\sF_t)_{t\geq 0}, \bP)$. 
 We call $(\Omega,\sF, (\sF_t)_{t\geq 0}, \bP; X,B)$ a weak solution of SDE \eqref{SDE} with starting point $x\in\mR^d$
if $B$ is an $\sF_t$-Brownian motion and
\begin{align}\label{Def9}
X_t=x+\int^t_0\sigma(X_s)\dif B_s+A^b_t,\ \forall t>0,\ \ \bP-a.s.,
\end{align}
where $A^b_t:=\lim_{n\to\infty}\int^t_0b_n(X_s)\dif s$ in the sense of u.c.p., 
and $b_n\in C^2(\mR^d)$ is any approximation sequence of $b$ so that for each $R>0$,
$$
\lim_{n\to\infty}\|(b_n-b)\chi_R\|_{-\alpha,p}=0.
$$
Here $A^b_t$ does not depend on the choice of approximation sequence $b_n\in C^2(\mR^d)$ of $b$.
\ed

\br\label{Re310}\rm
If $\bP\circ X^{-1}\in\sK^\alpha_p(\mC)$, then the above limit $A^b_t=\lim_{n\to\infty}\int^t_0b_n(X_s)\dif s$ does exist. 
Indeed, as in Proposition \ref{Pr13}, for any $f\in H^{-\alpha,p}_{loc}$,  there is a unique continuous 
$\sF_t$-adapted process $A^f_t$ such that for any $T, R>0$, 
$$
\bE\left(\sup_{t\in[0,T]}\left|\int^{t\wedge\eta_R}_0f_n(X_s)\dif s-A^f_t\right|^2\right)\leq C_{T,R}\|(f_n-f)\chi_R\|^2_{-\alpha,p},
$$
where $f_n:=f*\rho_n$ and $\eta_R:=\inf\{t>0:|X_t|>R\}$. Moreover, we also have for all $0\leq t_0<t_1\leq T$,
\begin{align}\label{Kr00}
\bE\left|A^f_{t_1\wedge\eta_R}-A^f_{t_0\wedge\eta_R}\right|^2\leq C_{T,R}(t_1-t_0)^{1+\gamma}\|f\chi_R\|^2_{-\alpha,p}.
\end{align}
\er

To show the equivalence between  weak solutions and martingale solutions, we  need the following It\^o's formula established by F\"ollmer in \cite{Fo}
and a stochastic version of Young's integral.
\bl\label{Ito}
Let $X_t=X_0+M_t+A_t$ be a Dirichlet process, where $M_t$ is a continuous local martingale, and $A_t$ is a locally zero energy process
(see Remark \ref{Re23} (i)).
For any $f\in C^2$, we have
$$
f(X_t)-f(X_0)= \int_0^t \p_i f(X_s) \dif M^i_s+ \int_0^t \p_i f(X_s) \dif A^i_s+ \frac{1}{2} \int_0^t \p_{ij}f(X_s) \dif [M^i, M^j]_s,
$$
where $ \int_0^t \p_i f(X_s) \dif A^i_s$ is defined as the limit in probability of the usual Riemmanian sum.
\el
Let $p\geq 1$ and $\beta>0$. For a stochastic process $A_t$ and $T>0$, we write
$$
\cH^{\beta,p}_T(A):=\|A_0\|_{L^p(\Omega)}+\sup_{s\not =t, s,t\in[0,T]}\frac{\|A_t-A_s\|_{L^p(\Omega)}}{|t-s|^{\beta}}.
$$
The following lemma is a slight extension of \cite[Lemma 2.2]{Ba-Ch1}, which can be considered as an analogue of usual Young's integral.
\bl\label{Le72}
Let $A,K$ be two stochastic processes and $p,q,r\in[1,\infty)$ with $\frac{1}{r}=\frac{1}{p}+\frac{1}{q}$. 
Suppose that for any $T>0$, there are $\gamma,\beta\in(0,1]$ with $\gamma+\beta>1$ such that
\begin{align*}
\cH^{\beta,p}_T(A)<\infty,\ \  \cH^{\gamma,q}_T(K)<\infty.
\end{align*}
For $n\in\mN$ and $s>0$, define $s_n:=[2^ns]2^{-n}$, where $[a]$ denotes the integer part of real number $a$.
Then for any $T>0$,
$\int^\cdot_0 K_{s_n}\dif A_s$ converges in $C([0,T];L^r(\Omega))$ and the limit is denoted by $\int^\cdot_0 K_s\dif A_s$. 
Moreover, there is a constant $C>0$
depending only on $\beta,\gamma$ and $T$ such that for all $t\in[0,T]$,
\begin{align}\label{Moment-J}
\left\|\int^t_0 K_{s_n}\dif A_s- \int^t_0 K_{s}\dif A_s\right\|_{L^r(\Omega)}\leq C\,\cH^{\gamma,q}_T(K)\,\cH^{\beta,p}_T(A)\,2^{-n(\beta+\gamma-1)},
\end{align}
and for all $0\leq t'<t\leq T$,
\begin{align}
\left\|\int^{t}_{t'} K_{s}\dif A_s\right\|_{L^r(\Omega)}\leq C\,\cH^{\gamma,q}_T(K)\,\cH^{\beta,p}_T(A)\,(t-t')^{\beta}.
\end{align}
\el 
\begin{proof}
For simplicity of notation, we write $\delta^n_k:=k2^{-n}$. Noticing that
$$
\int^t_0 K_{s_n}\dif A_s=\sum_{k=0}^{[2^nt]-1}  K_{\delta^n_k} \big(A_{\delta^n_{k+1}}-A_{\delta^n_k}\big)+
K_{[2^n t]2^{-n}} \big(A_{t}-A_{[2^n t]2^{-n}}\big),
$$
we have
\begin{align*}
\int^t_0(K_{s_{n+1}}-K_{s_n})\dif A_s=\sum_{k\hbox{\tiny{ even}}}^{2[2^{n}t]} \big(K_{\delta^{n+1}_{k+1}}-K_{\delta^{n+1}_k}\big) \big(A_{\delta^{n+1}_{k+2}}-A_{\delta^{n+1}_{k+1}}\big)+R_t^n,
\end{align*}
where 
\begin{align*}
R_t^n:=1_{\{2^nt-[2^nt]>1/2\}}\big(K_{[2^{n+1}t]2^{-n-1}} -K_{[2^nt]2^{-n}}\big)\big(A_t-A_{[2^{n+1}t]2^{-n-1}}\big).
\end{align*}
Suppose $2^{-m}<t-t'\leq 2^{1-m}$ for some $m\in\mN$. Then for $n\geq m$, we have
\begin{align*}
\int^{t}_{t'}(K_{s_{n+1}}-K_{s_n})\dif A_s=\sum_{k\ \mbox{\tiny even}, k>2[2^nt']}^{2[2^nt]} \big(K_{\delta^{n+1}_{k+1}}-K_{\delta^{n+1}_k}\big) 
\big(A_{\delta^{n+1}_{k+2}}-A_{\delta^{n+1}_{k+1}}\big)+R_{t}^n-R_{t'}^n,
\end{align*}
and by H\"older's inequality,
\begin{align*}
\left\|\int^{t}_{t'}(K_{s_{n+1}}-K_{s_n})\dif A_s\right\|_{L^r(\Omega)}\leq & 
\sum_{k\ \mbox{\tiny even}, k>2[2^nt']}^{2[2^nt]}  \|K_{\delta^{n+1}_{k+1}}-K_{\delta^{n+1}_k}\|_{L^q(\Omega)} 
\|A_{\delta^{n+1}_{k+2}}-A_{\delta^{n+1}_{k+1}}\|_{L^p(\Omega)}\\
&+ \big\| K_{[2^{n+1}t]2^{-n-1}} -K_{[2^nt]2^{-n}}\big\|_{L^q(\Omega)} \big\| A_t-A_{[2^{n+1}t]2^{-n-1}}\big\|_{L^p(\Omega)}\\
&+ \big\| K_{[2^{n+1}t']2^{-n-1}} -K_{[2^nt']2^{-n}}\big\|_{L^q(\Omega)} \big\| A_{t'}-A_{[2^{n+1}t']2^{-n-1}}\big\|_{L^p(\Omega)}\\
\leq & \cH^{\gamma,q}_T(K)\,\cH^{\beta,p}_T(A)\left(\sum_{k\ \mbox{\tiny even}, k>2[2^nt']}^{2[2^nt]}  2^{-n(\beta+\gamma)} + 2^{-n (\beta+\gamma)}\right)\\
\leq&\cH^{\gamma,q}_T(K)\,\cH^{\beta,p}_T(A)\, 3\cdot 2^{-n(\beta+\gamma-1)}(t-t').
\end{align*}
Moreover, since $2^{-m}<t-t'\leq 2^{1-m}$, it is easy to see that
$$
\left\|\int^t_{t'} K_{s_m}\dif A_s\right\|_{L^r(\Omega)}\leq 2\sup_{s\in[0,T]}\|K_s\|_{L^q(\Omega)}\cH^{\beta,p}_T(A)(t-s)^\beta.
$$
Combining the above two inequalities, we get the desired result.  
\end{proof}
Now we can show the following equivalence.
\bp\label{Pr29}
Let $\mP\in\sP(\mC)$ satisfy that for any $T,R>0$ and $s,t\in[0,T]$,
\begin{align}\label{IP1}
\mE |w_{t\wedge \tau_R}-w_{s\wedge \tau_R}|^2\leq C_{T,R}|t-s|.
\end{align}
Let $\alpha\in[0,1]$ and $p>1$. Assume that $b\in H^{-\alpha,p}_{loc}$ and $\sigma,\sigma^{-1}$ are locally bounded.
Then $\mP\in\sM^{\alpha,p}_{\sigma,b}(x)$
if and only if there is a weak solution $(\Omega,\sF, (\sF_t)_{t\geq 0}, \bP; X,B)$ in the sense of Definition \ref{Def99}
so that $\bP\circ X^{-1}=\mP\in\sK^\alpha_p(\mC)$.
\ep
\begin{proof}
(i) Let $(\Omega,\sF, (\sF_t)_{t\geq 0},\bP;X,B)$ be a weak solution of SDE \eqref{SDE} satisfying
\begin{align}\label{GD44}
\bP\circ X^{-1}\in\sK^\alpha_p(\mC).
\end{align}
By \eqref{Kr00}, $X$ is a Drichlet process. For any $f\in C^\infty$, by Lemma \ref{Ito} we have
$$
f(X_t)= f(x)+\int_0^t \nabla f(X_s)\cdot \sigma(X_s)\dif B_s+ \int_0^t \sL^\sigma f(X_s) \dif s+ \int_0^t \nabla f(X_s) \cdot \dif A_t^b, 
$$
where the last term in the right hand side is defined as the limit in probability 
of Riemannian sum (see Lemma \ref{Ito}).  
To show $\bP\circ X^{-1}\in \sM^{\alpha,p}_{\sigma,b}(x)$, by definition it suffices to prove that for any $t>0$,
\begin{align}\label{Chain}
\int_0^t \nabla f(X_s) \cdot \dif A_t^b=A^{b\cdot\nabla f}_t,\ \ \bP-a.s.,
\end{align}
where  the right hand side is defined as in Remark \ref{Re310}. 
By \eqref{GD44} and Remark \ref{Re310}, we have 
\begin{align}\label{NK1}
A^{b\cdot\nabla f}_t=\lim_{n\to \infty }\int_0^t\nabla f(X_s)\cdot \dif A^{b_n}_s
=\lim_{n\to \infty }\int_0^t(b_n\cdot\nabla f)(X_s)\dif s \mbox{  in probability}, 
\end{align}
where $b_n:=b*\rho_n$, and for any $T, R>0$ and $s,t\in[0,T]$,
\begin{align*}
\bE |A^{b_n}_{t\wedge \eta_R}-A^{b_n}_{s\wedge \eta_R}|^2
\leq C_{T,R} \|b_n\chi_R\|_{-\a, p}^2 |t-s|^{1+\gamma}\leq C'_{T,R} \|b\chi_R\|_{-\a, p}^2 |t-s|^{1+\gamma}, 
\end{align*}
where $\eta_R:=\inf\{t>0:|X_t|>R\}$, and also
\begin{align}\label{NL33}
\bE |A^{b}_{t\wedge \eta_R}-A^{b}_{s\wedge \eta_R}|^2\leq C'_{T,R} \|b\chi_R\|_{-\a, p}^2 |t-s|^{1+\gamma}.
\end{align}
By \eqref{Def9} and \eqref{NL33}, it is easy to see that for any $T,R>0$ and $s,t\in[0,T]$,
$$
\bE |\nabla f(X_{t\wedge \eta_R})-\nabla f(X_{s\wedge \eta_R})|^2 \leq \|\nabla^2f \|_\infty^2 \bE |X_{t\wedge \eta_R}-X_{s\wedge \eta_R}|^2 \leq C_{T,R} |t-s|.
$$
Hence, by \eqref{Moment-J} with $p=q=2$,
\begin{align*}
&\lim_{m\to\infty}\sup_{n\in\mN\cup\{\infty\}}\bE \left|\int_0^{t\wedge\eta_R} \nabla f(X_s)\cdot\dif A^{b_n}_s-\int^{t\wedge\eta_R}_0\nabla f(X_{s_m})\cdot\dif A^{b_n}_s\right|\\
&\quad\leq C\lim_{m\to\infty}\sup_{n\in\mN\cup\{\infty\}}\Big(\cH^{\frac{1}{2},2}_T(\nabla f(X_{\cdot\wedge\eta_R}))\,
\cH^{\frac{1+\gamma}{2},2}_T(A^{b_n}_{\cdot\wedge\eta_R})\,2^{-m\gamma/2}\Big)=0,
\end{align*}
where $b_\infty:=b$ and $s_m:=[2^m s]2^{-m}$. Since $\bP(\lim_{R\to\infty}\eta_R=\infty)=1$, we further have
\begin{align}\label{GR1}
\lim_{m\to\infty}\sup_{n\in\mN\cup\{\infty\}}\bP \left(\left|\int_0^t \nabla f(X_s)\cdot\dif A^{b_n}_s
-\int^t_0\nabla f(X_{s_m})\cdot\dif A^{b_n}_s\right|>\eps\right)=0,\ \forall \eps>0.
\end{align}
On the other hand, since $A^{b_n}\to A^b$ in the sense of u.c.p.,  for fixed $m\in\mN$,
by writing the integral as a discretization sum, we have
$$
\int_0^t \nabla f(X_{s_m})\cdot \dif A^{b_n}_s\stackrel{n\to\infty}{\longrightarrow}\int_0^t \nabla f(X_{s_m})\cdot \dif A^b_s\mbox{ in probability,}
$$
which, together with \eqref{GR1} and \eqref{NK1}, implies \eqref{Chain}. 
\\
\\
(ii) Suppose that $\mP\in\sM^{\alpha,p}_{\sigma,b}(x)$ satisfies \eqref{IP1}. 
By choosing $f(x)=x_i$ in \eqref{EQ1}, one sees that  $M^i_t:=w^i_{t}-x^i -A^{b^i}_{t}$ is a continuous local martingale under $\mP$. 
By Lemma \ref{Ito} again, we get 
$$
w_t^iw_t^j-x^ix^j= \int_0^t (w_s^j \dif M^i_s+w_s^i \dif M^j_s )+ \int_0^t (w_s^j \dif A^{b^i}_s+w_s^i \dif A_s^{b^j})+ [M^i, M^j]_t.
$$
On the other hand, for any $i,j=1,\cdots,d$, if we choose $f(x)=x_ix_j$ in \eqref{EQ1}, then
$$
w^i_{t}w^j_t-x^ix^j -A^{x^jb^i}_{t}-A^{x^ib^j}_{t}-\int_0^t a^{ij}(w_s)\dif s
$$
is also a continuous local martingale. As in showing \eqref{Chain}, by \eqref{IP1} and \eqref{Kr0}, we have
$$
A^{x^ib^j}_{t}=\int_0^tw_s^i \dif A_s^{b^j},\ \ i,j=1,\cdots,d.
$$
Hence, 
$$
[M^i, M^j]_t=\int_0^t a^{ij} (w_s)\dif s. 
$$
Now we define
$$
B_t:=\int^t_0\sigma^{-1}(w_s)\dif M_s,\ \ t\geq 0.
$$
Since $\sigma^{-1}$ is locally bounded,
$B$ is a continuous $\cB_t(\mC)$-local martingale under $\mP$ and by definition,
$$
[B^i, B^j]_t=\delta_{ij}\, t,\ i,j=1,\cdots, d.
$$
By L\'evy's characterization, $B$ is a $\cB_t(\mC)$-Brownian motion under $\mP$. Moreover,
$$
w_{t}=x+A^{b}_{t}+\int^t_0\sigma(w_s)\dif B_s,\ \ \mP-a.s.
$$
Thus $(\mC,\cB(\mC), (\cB_t(\mC))_{t\geq 0}, \mP; w,B)$ is a weak solution in the sense of Definition \ref{Def99}.
\end{proof}

\section{Cauchy problem for PDEs with distributional drifts}

In this section we solve the following Cauchy problem of PDEs with distributional drifts:
\begin{align}\label{PDE8}
\p_tu=\sL^a u-\lambda u+b\cdot\nabla u+f,\ \ u(0)=\varphi.
\end{align}
First of all we prepare two freezing lemmas in Bessel potential spaces for later use.
\bl\label{Unite}
Let $\phi$ be a nonnegative and nonzero smooth function with compact support.
Define $\phi_z(x):=\phi(x-z)$. For any $\alpha\in\mR$ and $p\in (1,\infty)$, there exists a constant $C\geq 1$ depending only on $\alpha,p,\phi$ such that 
for all $f\in H^{\alpha,p}$,
\be\label{GJ1}
C^{-1} \|f\|_{\alpha,p}\leq \left(\int_{\mR^d}  \|\phi_z f\|_{\alpha,p}^p\dif z\right)^{1/p} \leq C \|f\|_{\alpha,p}. 
\ee
\el
\bpf
Define 
$$
T^\phi f(z, x):=\phi_z(x)f(x),\ \ x,z\in\mR^d.
$$
Suppose that we have proved that for all $1<p<\infty$ and $\alpha\in\mR$, there is a $C>0$ such that
\begin{align}\label{EY1}
\|T^\phi f\|_{L^p(\mR^d; H^{\alpha,p})}\leq C\|f\|_{\alpha,p},
\end{align}
that is, the right hand side estimate in \eqref{GJ1} was proved,
then the left hand side estimate follows by a duality argument. In fact, letting $p':=\frac{p}{p-1}$, we have
\begin{align*}
\|f\|_{\alpha,p}&=\sup_{\|g\|_{-\a, p'}\leq 1}\left|\int_{\mR^d}f(x)\cdot g(x)\dif x\right|
=\sup_{\|g\|_{-\a, p'}\leq 1}\left|\int_{\mR^d}\!\int_{\mR^d}f(x)\cdot g(x)\phi^2_z(x)\dif z\dif x\right|\Big/\int_{\mR^d}\phi^2\\
&\lesssim \sup_{\|g\|_{-\a, p'}\leq 1}\left(\int_{\mR^d}\|\phi_z f\|_{\alpha,p}\|\phi_z g\|_{-\alpha,{p'}}\dif z\right)\leq\sup_{\|g\|_{-\a, p'}\leq 1}\|T^\phi f\|_{L^p(\mR^d; H^{\alpha,p})}\|T^\phi g\|_{L^{p'}(\mR^d; H^{-\alpha,{p'}})}\\
&\lesssim\|T^\phi f\|_{L^p(\mR^d; H^{\alpha,p})}\sup_{\|g\|_{-\a, p'}\leq 1}\|g\|_{-\alpha,{p'}}
=\|T^\phi f\|_{L^p(\mR^d; H^{\alpha,p})}.
\end{align*}

To show \eqref{EY1}, by a standard interpolation method, it suffices to prove it for $\alpha=0,\pm 2k,\cdots$. 
For $\alpha=0,2,4,\cdots$, it follows by the chain rule. 
For $\alpha=-2$, still by duality, we have
\begin{align*}
\|T^\phi f\|^p_{L^p(\mR^d; H^{-2,p})}&=\int_{\mR^d}\sup_{\|g\|_{2, p'}\leq 1}|\<\phi_zf, g\>|^p\dif z
=\int_{\mR^d}\sup_{\|g\|_{2, p'}\leq 1}|\<(\mI-\Delta)^{-1}f,(\mI-\Delta)(\phi_zg)\>|^p\dif z.
\end{align*}
Recalling that
$$
(\mI-\Delta)(\phi_z g)=-\Delta\phi_z\cdot g+\phi_z\cdot (\mI-\Delta) g-2\nabla\phi_z\cdot\nabla g,
$$
we have
\begin{align*}
&|\<(\mI-\Delta)^{-1}f,(\mI-\Delta)(\phi_z g)\>|
\leq\|(\mI-\Delta)^{-1}f\cdot\Delta\phi_z\|_p\|g\|_{p'}\\
&\qquad+\|(\mI-\Delta)^{-1}f\cdot\phi_z\|_p\|(\mI-\Delta)g\|_{p'}+2 \|(\mI-\Delta)^{-1}f\cdot \nabla\phi_z\|_p\|\nabla g \|_{p'}\\
&\quad\lesssim\Big(\|(\mI-\Delta)^{-1}f\cdot\Delta\phi_z\|_p
+\|(\mI-\Delta)^{-1}f\cdot\phi_z\|_p+\|(\mI-\Delta)^{-1}f\cdot \nabla\phi_z\|_p\Big)\|g\|_{2, p'}.
\end{align*}
Combining the above inequalities, we get 
\begin{align*}
\|T^\phi f\|^p_{L^p(\mR^d; H^{-2,p})}
&\lesssim\int_{\mR^d}\Big(\|(\mI-\Delta)^{-1}f\cdot\Delta\phi_z\|_p
+\|(\mI-\Delta)^{-1}f\cdot\phi_z\|_p+\|(\mI-\Delta)^{-1}f\cdot  \nabla\phi_z \|_p\Big)^p\dif z\\
&\lesssim\|(\mI-\Delta)^{-1}f\|^p_p=\|f\|^p_{-2,p}.
\end{align*}
For general $\alpha=-4,-6,\cdots$, it follows by similar calculations.
\epf
\br\rm
In fact, a discretized version of Lemma \ref{Unite} was proven in \cite{Kr1}. Our proof presented here is much simpler. 
\er
\bl
For any  $\gamma\in(0,1)$, $p\in(d/\gamma,\infty)$ and $\alpha\in[0,\gamma]$, 
there is a constant $C=C(d,\alpha,\gamma,p)>0$ such that for all $\delta\in(0,1)$ and $z\in\mR^d$,
\begin{align}\label{RG3}
\|(f_z-f)\chi^\delta_z\|_{\alpha,p}\leq C\|\Delta^{\gamma/2}f\|_{p}\delta^{\gamma-\alpha},
\end{align}
where $f_z:=f(z)$ and $\chi^\delta_z(x):=\chi_\delta(x-z)$, $\chi_\delta$ is defined as \eqref{Cut}. 
Moreover, for any continuous function $f$ with $\|\Delta^{\gamma/2}f\|_p<\infty$, we have
\begin{align}\label{RG2}
\lim_{\delta\to 0}\sup_z\|(f_z-f)\chi^\delta_z\|_{\alpha,p}=0.
\end{align}
\el
\bpf
(i) We first consider the case $\alpha<\gamma$. 
Since $p\gamma>d$, by \eqref{RG1} we have
\begin{align}\label{YT2}
\|(f_z-f)\chi^\delta_z\|^{p}_{p}\lesssim\|\Delta^{\gamma/2}f\|^p_{p}\delta^{\gamma p-d}\int_{\mR^d}|\chi_\delta(x-z)|^{p}\dif x
\lesssim\| \Delta^{\gamma/2} f\|^p_{p}\delta^{\gamma p}.
\end{align}
Noticing that 
$$
\|\Delta^{\gamma/2}\chi^\delta_z\|_{p}=\|\Delta^{\gamma/2}\chi_\delta\|_{p}\leq\delta^{-\gamma}\|(\Delta^{\gamma/2}\chi)(\cdot/\delta)\|_{p}\lesssim\delta^{-\gamma+d/p},
$$
by \eqref{GF1} and \eqref{RG1} we have for $|y|\leq\delta$,
\begin{align*}
\|((f_z-f)\chi^\delta_z)(\cdot+y)-(f_z-f)\chi^\delta_z\|_{p}&\lesssim\|f(\cdot+y)-f\|_{p}\|\chi^\delta_z\|_\infty+\|(f_z-f)(\cdot+y)(\chi^\delta_z(\cdot+y)-\chi^\delta_z)\|_{p}\\
&\lesssim|y|^\gamma\|\Delta^{\gamma/2}f\|_{p}+\|\Delta^{\gamma/2}f\|_{p}\delta^{\gamma-d/p}\|\chi^\delta_z(\cdot+y)-\chi^\delta_z\|_{p}\\
&\lesssim|y|^\gamma\|\Delta^{\gamma/2}f\|_{p} (1+\delta^{\gamma-d/{p}}\|\Delta^{\gamma/2}\chi^\delta_z\|_{p})\lesssim |y|^\gamma \|\Delta^{\gamma/2}f\|_{p}.
\end{align*}
Thus, by \eqref{Fr} and \eqref{YT2} we have
\begin{align*}
\|\Delta^{\alpha/2}((f_z-f)\chi^\delta_z)\|_{p}&\leq\int_{\mR^d}\|(f_z-f)\chi^\delta_z)(\cdot+y)-(f_z-f)\chi^\delta_z\|_{p}\frac{\dif y}{|y|^{d+\alpha}}\\
&\lesssim  \|\Delta^{\gamma/2}f\|_{p}\int_{|y|\leq\delta}|y|^\gamma\frac{\dif y}{|y|^{d+\alpha}}
+ \|\Delta^{\gamma/2}f\|_{p}\delta^{\gamma}\int_{|y|>\delta}\frac{\dif y}{|y|^{d+\alpha}}\\
&\lesssim \delta^{\gamma-\alpha} \|\Delta^{\gamma/2}f\|_{p}.
\end{align*}
(ii) Next we consider the case $\alpha=\gamma$. By definition, we have
\begin{align*}
\Delta^{\gamma/2}((f_z-f)\chi^\delta_z)(x)=\Delta^{\gamma/2}f(x)\cdot\chi^\delta_z(x)
+\int_{\mR^d}(f(z)-f(x+y))(\chi^\delta_z(x+y)-\chi^\delta_z(x))\frac{\dif y}{|y|^{d+\gamma}}.
\end{align*}
Clearly, we have
$$
\|\Delta^{\gamma/2}f\cdot\chi^\delta_z\|_p\lesssim\|\Delta^{\gamma/2}f\|_{p}.
$$
To estimate the second term denoted by $\sI^\delta_z(x)$, noticing that for $|y|\leq\delta$,
$$
\chi^\delta_z(x+y)-\chi^\delta_z(x)=0\ \ \mbox{if}\ |x-z|>3\delta,
$$
we may write
\begin{align*}
\|\sI^\delta_z\|_p^p&=\int_{|x-z|\leq 3\delta}\left|\int_{|y|\leq\delta}(f(z)-f(x+y))(\chi^\delta_z(x+y)-\chi^\delta_z(x))\frac{\dif y}{|y|^{d+\gamma}}\right|^p\dif x\\
&+\int_{|x-z|\leq 3\delta}\left|\int_{|y|>\delta}(f(z)-f(x+y))(\chi^\delta_z(x+y)-\chi^\delta_z(x))\frac{\dif y}{|y|^{d+\gamma}}\right|^p\dif x\\
&+\int_{|x-z|>3\delta}\left|\int_{|y|>\delta}(f(z)-f(x+y))(\chi^\delta_z(x+y)-\chi^\delta_z(x))\frac{\dif y}{|y|^{d+\gamma}}\right|^p\dif x\\
&=:I_1+I_2+I_3.
\end{align*}
For $I_1$, by \eqref{RG1} we have
\begin{align*}
I_1&\lesssim\|\Delta^{\gamma/2}f\|^p_{p}\delta^{\gamma p-d}\int_{|x-z|\leq 3\delta}\left|\int_{|y|\leq\delta}|\chi^\delta_z(x+y)-\chi^\delta_z(x)|\frac{\dif y}{|y|^{d+\gamma}}\right|^p\dif x\\
&\lesssim\|\Delta^{\gamma/2}f\|^p_{p}\delta^{\gamma p-d}\int_{|x-z|\leq 3\delta}\left|\delta^{-1}\int_{|y|\leq\delta}\frac{\dif y}{|y|^{d+\gamma-1}}\right|^p\dif x\lesssim\|\Delta^{\gamma/2}f\|^p_{p},
\end{align*}
and
\begin{align*}
I_{2}\lesssim&  \|\Delta^{\gamma/2}f\|^p_{p}\int_{|x-z|\leq 3\delta}  \left(\int_{|y|>\delta} |x+y-z|^{\gamma-d/p}  
\frac{\dif y}{|y|^{d+\gamma} }\right)^p  \dif x\\
\lesssim& \|\Delta^{\gamma/2}f\|_{ p}^p \delta^d\  \left(\int_{|y|>\delta}|y|^{-d-d/p} \dif y\right)^p
\lesssim\|\Delta^{\gamma/2}f\|_{ p}^p. 
\end{align*}
For $I_3$, noticing that if $|x-z|>3\delta$, then $\chi_z^\delta(x)=0$ and if $|x+y-z|\leq 2\delta$, then $\chi_z^\delta(x+y)=0$, we have 
\begin{align*}
I_3&=\int_{|x-z|>3\delta}\left|\int_{|y|>\delta, |x+y-z|\leq 2\delta}(f(z)-f(x+y))\chi^\delta_z(x+y)\frac{\dif y}{|y|^{d+\gamma}}\right|^p\dif x\\
&\lesssim \int_{|x-z|>3\delta} \left| \int_{|y|\leq 2\delta} \|\Delta^{\gamma/2}f\|_{p} 
\delta^{\gamma-d/p}\frac{ \dif y}{(|x-z|-2\delta)^{d+\gamma}}\right|^p  \dif x
\lesssim \|\Delta^{\gamma/2}f\|_{p}^p.
\end{align*}
Combining the above calculations, we obtain \eqref{RG3}.
\\
\\
(iii) If $\alpha<\gamma$,  the limit \eqref{RG2} is obvious. If $\alpha=\gamma$, letting $f^n=f*\rho_n$, by \eqref{RG3} we have
\begin{align*}
\|(f_z-f)\chi^\delta_z\|_{\alpha,p}&\leq\|(f^n_z-f^n)\chi^\delta_z\|_{\alpha,p}+\|(f-f^n)_z-(f-f^n))\chi^\delta_z\|_{\alpha,p}\\
&\leq C\left(\|\nabla f^n\|_{p}\delta^{1-\alpha}+\|\Delta^{\alpha/2}(f-f^n)\|_{p}\right),
\end{align*}
which gives \eqref{RG2} by first letting $\delta\to 0$, then $n\to\infty$.
\epf

For $T>0$ and $\alpha\in\mR, 1<p<\infty$, we introduce the following Banach space:
$$
\mH^{\alpha,p}_T:=L^p([0,T]; H^{\alpha,p}).
$$
We first show the following result about constant coefficient equation.
\bt
Let $a(x)=a$ be a constant symmetric positive definite matrix. Let $\lambda>0$, $\alpha\in\mR$, $p>1$ and $T>0$.
For any $\varphi\in H^{2+\alpha,p}$ and $f\in \mH^{\alpha,p}_T$, there is a unique solution $u\in \mH^{2+\alpha,p}_T$ to the Cauchy problem
\begin{align}\label{KL1}
\p_t u=\sL^a u-\lambda u+f,\ \ u(0)=\varphi.
\end{align}
Moreover, for any $\theta\in[0,2]$, there is a constant $C>0$ only depending on the elliptic constant of $a$ 
and $\theta,p,d, T$ such that for all $\lambda\geq 1$,
\begin{align}\label{YT1}
\lambda^{1-\frac{\theta}{2}}\|u\|_{\mH^{\theta+\alpha,p}_T}\leq C\Big(\lambda^{1-\frac{\theta}{2}-\frac{1}{p}}
\|\varphi\|_{\theta+\alpha,p}+\|f\|_{\mH^{\alpha,p}_T}\Big).
\end{align}
\et
\begin{proof}
Let $P^a_t f(x):=\mE f(\sqrt{2a}\cdot B_t+x)$ be the Gaussian heat semigroup with diffusion matrix $\sqrt{2a}$.
By Duhamel's formula, the unique solution of \eqref{KL1} can be written as 
$$
u(t,x)=\e^{-\lambda t}P^a_t\varphi(x)+\int^t_0\e^{-\lambda(t-s)} P^a_{t-s}f(s,x)\dif s.
$$
By \cite[Theorem 1.1]{Kr0}, we have
\begin{align*}
\|\nabla^2 u\|^p_{\mH^{\alpha,p}_T}&\lesssim \int^T_0\|\nabla^2P^a_t\varphi\|^p_{\alpha,p}\dif t
+\int^T_0\left\|\nabla^2\int^t_0 P^a_{t-s}f(s,\cdot)\dif s\right\|^p_{\alpha,p}\dif t\lesssim\|\varphi\|^p_{2+\alpha,p}+\|f\|^p_{\mH^{\alpha,p}_T}.
\end{align*}
On the other hand, for any $\theta\in [0,2)$, noticing that
$$
\|P^a_tf\|_{\theta,p}\leq Ct^{-\theta/2}\|f\|_p,\ t>0,
$$
we have
\begin{align*}
\|u\|^p_{\mH^{\theta+\alpha,p}_T}&\lesssim\int^T_0\e^{-\lambda p t}\|P^a_t\varphi\|^p_{\theta+\alpha,p}\dif t
+\int^T_0\left(\int^t_0\e^{-\lambda(t-s)}\|P^a_{t-s}f(s,\cdot)\|_{\theta+\alpha,p}\dif s\right)^p\dif t\\
&\lesssim\left(\int^T_0\e^{-\lambda p t }\dif t\right) \|\varphi\|_{\theta+\alpha,p}^p
+\int^T_0\left(\int^t_0\e^{-\lambda(t-s)}(t-s)^{-\theta/2}\|f(s,\cdot)\|_{\alpha,p}\dif s\right)^p\dif t\\
&=\frac{1-\e^{-\lambda pT}}{\lambda p}\|\varphi\|_{\theta+\alpha,p}^p
+\int^T_0\left(\int^T_0\e^{-\lambda s}s^{-\theta/2}\|f(t-s,\cdot)\|_{\alpha,p}1_{\{t-s>0\}}\dif s\right)^p\dif t\\
&\lesssim\lambda^{-1}\|\varphi\|^p_{\theta+\alpha,p}+\lambda^{\frac{\theta p}{2}-p}\|f\|^p_{\mH^{\alpha,p}_T},
\end{align*}
where the last step is due to Minkowskii's inequality.
The proof is complete.
\end{proof}

To show the corresponding result for variable coefficient $a=\sigma\sigma^*/2$,
we make the following assumptions about $\sigma$:  
\begin{enumerate}[\bf (H$^\sigma_{\beta,q}$)]
\item $\|\Delta^{\beta/2}\sigma\|_{q}<\infty$ for some $\beta\in(0,1]$ and $q\in(\frac{d}{\beta},\infty)$,
and there is a constant $c_0\geq 1$ such that
\begin{align}\label{DG1}
c^{-1}_0|\xi|^2\leq|\sigma(x)\xi|^2\leq c_0|\xi|^2,\ \forall x,\xi\in\mR^d.
\end{align}
\end{enumerate}

\bt\label{Th43}
Let $\beta\in(0,1)$ and $q\in(\frac{d}{\beta},\infty)$.
Under {\bf (H$^\sigma_{\beta,q}$)}, for any $\alpha\in[0,\beta]$ and $p>\frac{d}{d-\a}$, 
there is a $\lambda_0\geq 1$ large enough such that 
for all $\lambda\geq\lambda_0$, $T>0$ 
and any $\varphi\in H^{2-\alpha,p}$, $f\in \mH^{-\alpha,p}_T$, 
there is a unique $u\in \mH^{2-\alpha,p}_T$ solving the following PDE
\begin{align}\label{PDE0}
u(t)=\varphi+\int^t_0\Big[(\sL^a -\lambda) u(s)+f(s)\Big]\dif s\ \mbox{ in  $H^{-\alpha,p}$}.
\end{align}
Moreover, for any $\theta\in[0,2]$,
there is a constant $C>0$ such that for all $\lambda\geq \lambda_0$,
\begin{align}\label{ER8}
\lambda^{1-\frac{\theta}{2}}\|u\|_{\mH^{\theta-\alpha,p}_T}\leq C\Big(\lambda^{1-\frac{\theta}{2}-\frac{1}{p}}\|\varphi\|_{2-\alpha,p}+\|f\|_{\mH^{\alpha,p}_T}\Big).
\end{align}
\et
\bpf
By the  standard continuity method, it suffices to show the a priori estimate \eqref{ER8}. We use the freezing coefficient argument.
Let $\phi$ be a nonnegative and nonzero smooth function with support in 
$B_1:=\{x\in\mR^d: |x|\leq 1\}$ 
and define for $z\in\mR^d$,
$$
\phi^\delta(x):=\delta^{-d}\phi(x\delta^{-1}),\ \phi^\delta_z(x):=\phi^\delta(x-z),\ \ a_z:=a(z).
$$ 
Multiplying both sides of PDE \eqref{PDE0} by $\phi^\delta_z$, we have
\begin{align*}
&\p_t(\phi^\delta_z u)= \sL^{a_z}(\phi^\delta_z u)-\lambda\phi^\delta_z u+
f\phi^\delta_z+(\sL^{a_z}-\sL^a) (\phi^\delta_zu)-\sL^a\phi^\delta_z\cdot u-2a^{ij}\p_i\phi^\delta_z \p_ju,
\end{align*}
where the above equality holds in $H^{-\alpha,p}$ for Lebesgue almost all $t\geq 0$.
Let $\chi^\delta_z(x):=\chi_\delta(x-z)$, where $\chi_\delta$ is defined as in \eqref{Cut}. Noticing $q>\frac{d}{\beta}$, 
we can choose $\gamma\in [\a, \beta]$ and $p_2\geq q$ such that 
$$
\tfrac{\gamma }{q\beta}=\tfrac{1}{p_2}<\tfrac{\a}{d}<1-\tfrac{1}{p}. 
$$
Since $\chi^\delta_z\equiv 1$ on the support of $\phi^\delta_z$, we have $\chi^\delta_z\cdot\nabla^m\phi^\delta_z=\nabla^m\phi^\delta_z$ for $m=0,1,2$. Thus, for any $\theta\in[0,2]$, by \eqref{YT1} and Lemma \ref{Prod},
\begin{align*}
\lambda^{1-\frac{\theta}{2}}\|\phi^\delta_z u\|_{\mH^{\theta-\alpha,p}_T}
&\lesssim \lambda^{1-\frac{\theta}{2}-\frac{1}{p}}\|\varphi\phi^\delta_z\|_{\theta-\alpha,p}
+\|f\phi^\delta_z\|_{\mH^{-\alpha,p}_T}+\| (\sL^{a_z}-\sL^a) (\phi^\delta_zu)\|_{\mH^{-\alpha,p}_T}\\
&\quad+\|\sL^a\phi^\delta_z\cdot u\|_{\mH^{-\alpha,p}_T}+\|a^{ij}\p_i\phi^\delta_z \p_ju\|_{\mH^{-\alpha,p}_T}\\
&\lesssim \lambda^{1-\frac{\theta}{2}-\frac{1}{p}}\|\varphi\phi^\delta_z\|_{\theta-\alpha,p}+
\|f\phi^\delta_z\|_{\mH^{-\alpha,p}_T}+\|(a_z-a)\chi^\delta_z\|_{\alpha,p_2} \|\nabla^2(\phi^\delta_zu)\|_{\mH^{-\alpha,p}_T}\\
&\quad+\|a^{ij}\chi^\delta_z\|_{\alpha,p_2} \Big(\|\p_{ij}\phi^\delta_z\cdot u\|_{\mH^{-\alpha,p}_T}+\|\p_i\phi^\delta_z\p_j u\|_{\mH^{-\alpha,p}_T}\Big).
\end{align*}
Here and below, the constant contained in $\lesssim$ is independent of $\lambda,\delta$ and $\eps$.
Since $\sigma$ is bounded and $\|\Delta^{\beta/2}\sigma\|_{q}<\infty$, by \eqref{Fr}, \eqref{GF1} and \eqref{RG1}, one sees that
$$
\|a\|_\infty+ \|\Delta^{\beta/2}a\|_{q}<\infty,
$$
and by Gagliardo-Nirenberg's inequality \eqref{Ga},
$$
\|\Delta^{\gamma/2}a\|_{p_2} \lesssim \|a\|_\infty^{1-\gamma/\beta} \|\Delta^{\beta/2} a \|_{q}^{\gamma/\beta}<\infty.
$$
Thus for any $\eps>0$, by \eqref{RG2}, we can choose $\delta$ small enough so that 
$$
\sup_z\|(a_z-a)\chi^\delta_z\|_{\alpha,p_2}\leq\eps.
$$
Moreover, it is easy to see that
$$
C_\delta:=\sup_z \|a\chi^\delta_z\|_{\alpha,p_2} <\infty.
$$
Combining the above calculations, we get that for any $\eps>0$, there is a $\delta>0$ such that
\begin{align*}
\lambda^{1-\frac{\theta}{2}}\|\phi^\delta_z u\|_{\mH^{\theta-\alpha,p}_T}
&\lesssim\lambda^{1-\frac{\theta}{2}-\frac{1}{p}}\|\varphi\phi^\delta_z\|_{\theta-\alpha,p}+\|f\phi^\delta_z\|_{\mH^{-\alpha,p}_T}+\eps\|\phi^\delta_zu\|_{\mH^{2-\alpha,p}_T}\\
&\quad+C_\delta\sum_{ij}\Big(\|\p_{ij}\phi^\delta_z\cdot u\|_{\mH^{-\alpha,p}_T}+\|\p_i\phi^\delta_z\p_j u\|_{\mH^{-\alpha,p}_T}\Big).
\end{align*}
Taking $p$-order power for both sides and then integrating with respect to $z$ and by Lemma \ref{Unite}, we get
$$
\lambda^{1-\frac{\theta}{2}} \|u\|_{\mH^{\theta-\alpha,p}_T}\lesssim\lambda^{1-\frac{\theta}{2}-\frac{1}{p}}\|\varphi\|_{\theta-\alpha,p}+
 \|f\|_{\mH^{-\alpha,p}_T}+\eps \|u\|_{\mH^{2-\alpha,p}_T}+C_\delta\|u\|_{\mH^{1-\alpha,p}_T}.
$$
Letting $\theta=2, 1$, respectively, 
and choosing first $\eps$ small enough and then $\lambda$ large enough, we obtain the desired estimate.
\epf

As an easy corollary of the above result, 
we have
\bt\label{Th404}
Let $\alpha\in(0,\frac{1}{2}]$, $p>\frac{d}{1-\alpha}$ and $\beta\in[\alpha,1]$, $q\in(\frac{d}{\beta},\infty)$. 
Under  {\bf (H$^\sigma_{\beta,q}$)} and $b\in H^{-\a, p}$, there is a $\lambda_0\geq 1$ large enough such that 
for all $\lambda\geq\lambda_0$, $T>0$ and 
 any $\varphi\in H^{2-\alpha,p}$, $f\in \mH^{-\alpha,p}_T$, 
there is a unique $u\in \mH^{2-\alpha,p}_T$  such that 
\begin{align}\label{PDE00}
u(t)=\varphi+\int^t_0\Big[(\sL^a -\lambda+b\cdot\nabla) u(s)+f(s)\Big]\dif s\ \mbox{ in  $H^{-\alpha,p}$}.
\end{align}
Moreover, for any $\theta\in[0,2]$, there is a constant $C>0$ 
which only depends on the parameters and the constants in the assumptions,
\begin{align}\label{ER888}
\lambda^{1-\frac{\theta}{2}}\|u\|_{\mH^{\theta-\alpha,p}_T}\leq C\Big(\lambda^{1-\frac{\theta}{2}-\frac{1}{p}}\|\varphi\|_{2-\alpha,p}+\|f\|_{\mH^{-\alpha,p}_T}\Big).
\end{align}
\et
\bpf
By the continuity method, we still only need to prove \eqref{ER888}. Let $b_n =b*\rho_n$. 
By \eqref{ER8} and \eqref{EW2} with $p_1=p$ and $p_2=pd/(d-p(1-2\alpha))>d/\alpha$, we have
\begin{align*}
&\lambda^{1-\frac{\theta}{2}}\|u\|_{\mH^{\theta-\alpha,p}_T}\leq C\left(\lambda^{1-\frac{\theta}{2}-\frac{1}{p}}\|\varphi\|_{2-\alpha,p}
+\|f+b\cdot\nabla u\|_{\mH^{-\alpha,p}_T}\right)\\
&\quad\leq C\lambda^{1-\frac{\theta}{2}-\frac{1}{p}}\|\varphi\|_{2-\alpha,p}
+\|f\|_{\mH^{-\alpha,p}_T}+C\|b-b_n\|_{-\alpha,p}\|\nabla u\|_{\mH^{\alpha,p_2}_T}+C\|b_n\cdot\nabla u\|_{\mL^{p}_T}\\
&\quad\stackrel{\eqref{Fr2}}{\leq} C\lambda^{1-\frac{\theta}{2}-\frac{1}{p}}\|\varphi\|_{2-\alpha,p}
+\|f\|_{\mH^{-\alpha,p}_T}+C\|b-b_n\|_{-\alpha,p}\|u\|_{\mH^{2-\alpha,p}_T}+C\|b_n\|_{\infty}\|u\|_{\mH^{1,p}_T}.
\end{align*}
First choosing $\theta=2$ and $n$ large enough so that $C\|b-b_n\|_{-\alpha,p}\leq 1/2$, then letting $\theta=1+\alpha$
and $\lambda$ large enough so that $C\|b_n\|_{\infty}\leq \lambda^{\frac{1-\alpha}{2}}/4$, we get \eqref{ER888}.
\epf
\br\label{Re47}\rm
Notice that $u$ satisfies \eqref{PDE00} if and only if $u_\lambda(t,x):=\e^{\lambda t}u(t,x)$ satisfies
$$
u_\lambda(t)=\varphi+\int^t_0\Big[(\sL^a +b\cdot\nabla)u_\lambda(s)+\e^{\lambda s}f(s)\Big]\dif s\ \mbox{ in  $H^{-\alpha,p}$}.
$$
\er

We also have the following solvability to elliptic equations.
\bt\label{Th44}
Assume that one of the following two conditions holds:
\begin{enumerate}[(i)]
\item Let $\beta\in(0,1)$, $q\in(\frac{d}{\beta},\infty)$ and $\alpha\in(0,\beta]$, $p>\tfrac{d}{d-\a}$.
Assume that $b=0$ and {\bf (H$^\sigma_{\beta,q}$)} hold. 
\item Let  $\alpha\in(0,\frac{1}{2}]$, $p>\frac{d}{1-\alpha}$ and $\beta\in[\alpha,1]$, $q\in(\frac{d}{\beta},\infty)$.
Assume that $b\in H^{-\a, p}$ and  {\bf (H$^\sigma_{\beta,q}$)} hold.
\end{enumerate}
Then there is a $\lambda_0\geq 1$ large enough such that for all $\lambda\geq\lambda_0$ and any $f\in H^{-\alpha,p}$, 
there is a unique $u\in H^{2-\alpha,p}$ 
such that 
\begin{align}\label{PDE}
(\sL^a-\lambda+b\cdot\nabla)u=f\ \ \mbox{ in $H^{-\alpha,p}$}.
\end{align}
Moreover, for any $\theta\in[0,2]$,
there is a constant $C>0$ which only depends on the parameters and the constants in the assumptions such that for all $\lambda\geq \lambda_0$,
\begin{align}\label{ER88}
\lambda^{1-\frac{\theta}{2}} \|u\|_{\theta-\alpha,p}\leq C\|f\|_{-\alpha,p}.
\end{align}
\et
\bpf
We only consider the case (ii). Case (i) is similar by Theorem \ref{Th43}.
Let $T>0$ and 
$\phi:\mR\to\mR$ be a nonzero smooth function with compact support in $(0,T)$. 
Let $u\in H^{2-\alpha,p}$ solve elliptic equation \eqref{PDE}.
Then $\bar u(t,x):=u(x)\phi(t)$ satisfies the parabolic equation
$$
\p_t\bar u=(\sL^a -\lambda+b\cdot\nabla) \bar u+ f\phi+u\phi'\ \mbox{ in  $H^{-\alpha,p}$}.
$$
For any $\theta\in[0,2]$, by \eqref{ER888} we have 
$$
\lambda^{1-\frac{\theta}{2}}\|\bar u\|_{\mH^{\theta-\alpha,p}_T}\leq C\|f\phi+u\phi'\|_{\mH^{\alpha,p}_T}.
$$
Hence,
$$
\lambda^{1-\frac{\theta}{2}}\|u\|_{\theta-\alpha,p}\|\phi\|_{L^p(0,T)}\leq C\Big(\|f\|_{-\alpha,p}\|\phi\|_{L^p(0,T)}+\|u\|_{-\alpha,p}\|\phi'\|_{L^p(0,T)}\Big)
$$
and
$$
\lambda^{1-\frac{\theta}{2}}\|u\|_{\theta-\alpha,p}\leq C_1\|f\|_{-\alpha,p}+C_2\|u\|_{-\alpha,p}.
$$
Choosing $\theta=0$ and $\lambda\geq 2C_2$,
we get the desired estimate.
\epf

\section{Main results and proofs}
\subsection{Statement of main results}

We make the following assumptions about $b$: 
\begin{enumerate}[\bf (H$^b_{\alpha,p}$)]
\item $b=b^{(1)}+b^{(2)}$, where $b^{(1)}$ satisfies that for some $\vartheta\geq 0$ and $\kappa_0,\kappa_1,\kappa_2>0$,
\begin{align}\label{Diss} 
\frac{\<x,b^{(1)}(x)\>}{\sqrt{1+|x|^2}}\leq -\kappa_0|x|^{\vartheta}+\kappa_1,\ \ |b^{(1)}(x)|\leq \kappa_2(1+|x|^\vartheta),
\end{align}
and $b^{(2)}\in H^{-\alpha,p}$ for some $\alpha\in(0,\frac{1}{2}]$ and $p\in(\frac{d}{1-\alpha},\infty)$.
\end{enumerate}

Our first main result is
\bt\label{TH0}
Let $\alpha\in(0,\frac{1}{2}]$, 
$p\in(\frac{d}{1-\alpha},\infty)$ and $\beta\in[\alpha,1]$, $q\in(\frac{d}{\beta},\infty)$. Under {\bf (H$^\sigma_{\beta,q}$)} and {\bf (H$^b_{\alpha,p}$)},
for any $x\in\mR^d$, there exists a unique martingale solution $\mP_x\in\sM^{\alpha,p}_{\sigma,b}(x)$ 
to SDE \eqref{SDE}. Moreover, letting $\mE_x:=\mE^{\mP_x}$, we have the following conclusions:
\begin{enumerate}[(i)]
\item For any $T>0$ and $m\in\mN$, there is a constant $C_{T}>0$ such that for all $0\leq t_0<t_1\leq T$,
\begin{align}
\mE_x|w_{t_1}-w_{t_0}|^{2m}\leq C_{T}(t_1-t_0)^m,
\end{align}
and for all $f\in H^{-\alpha,p}$,
\begin{align}\label{GD1}
\mE_x\left|A^f_{t_1}-A^f_{t_0}\right|^{2m}\leq C_{T}(t_1-t_0)^{(2-\alpha-\frac{d}{p})m}\|f\|^{2m}_{-\alpha,p}.
\end{align}

\item
 If $\vartheta=0$ in \eqref{Diss}, then for any $\varphi\in H^{2-\a, p}$, $u(t,x):=P_t\varphi(x):=\mE_x\varphi(w_t)\in L^p_{loc}(\mR_+; H^{2-\alpha,p})$ 
uniquely solves the following Cauchy problem in $H^{-\alpha,p}$,
\begin{align}\label{PDE9}
\p_t u=(\sL^a+b\cdot\nabla)u,\ \ u(0)=\varphi. 
\end{align}
Moreover, $P_t$ admits a density $p(t,x,y)$ enjoying the following two-sided estimate: for some $c_1,c_2\geq 1$
and all $t>0$, $x,y\in\mR^d$,
\begin{align}\label{Es5}
c_1^{-1}t^{-d/2}\e^{-c_2 |x-y|^2/t}\leq p(t,x,y)\leq c_1t^{-d/2}\e^{-c_2^{-1} |x-y|^2/t},
\end{align}
and gradient estimate:  for some $c_3,c_4>0$ and all $t>0$, $x,y\in\mR^d$,
\begin{align}\label{Es6}
|\nabla_x p(t,x,y)|\leq c_3t^{-(d+1)/2}\e^{-c_4 |x-y|^2/t}.
\end{align}
\item If $\vartheta>0$ in \eqref{Diss}, then $P_t$ admits a unique invariant probability measure $\mu(\dif x)=\varrho(x)\dif x$ with $\varrho\in H^{\gamma,r}$, where $\gamma\in(0,\beta\wedge(1-\alpha)]$ and $r\in(1,\frac{d}{d+\gamma-1})$.
\end{enumerate}
\et

\br\rm
The above (iii) seems to be new even for $b^{(2)}\equiv 0$ although there are systematic studies
about the regularity of invariant measures $\varrho$  in the monograph \cite{Bo-Kr-Ro-Sa}. 
\er

As an easy corollary of Theorem \ref{TH0} and Proposition \ref{Pr29}, we have
\bc
Under the same assumptions of Theorem \ref{TH0}, there exists a unique weak solution 
$(\Omega,\sF, (\sF_t)_{t\geq 0}, \bP; X,B)$ for SDE \eqref{SDE} so that $\bP\circ X^{-1}\in\sK^\alpha_p(\mC)$.
\ec

In the above corollary, we require that the law of weak solution satisfies the local Krylov estimate, that is, $\bP\circ X^{-1}\in\sK^\alpha_p(\mC)$. 
This is crucial when we use Zvonkin's transformation to show the uniqueness. 
Nevertheless, under some extra assumptions (see \eqref{Def88} below), 
we can directly prove such a priori estimate for any weak solutions as stated below.
\bt\label{TH29}
Let $\alpha\in(0,\frac{1}{2})$, $p\in(\frac{d}{1/2-\alpha},\infty)$ and $\beta\in[\alpha,1]$, $q\in(\frac{d}{\beta},\infty)$.
Under {\bf (H$^\sigma_{\beta,q}$)} and {\bf (H$^b_{\alpha,p}$)}, for any $x\in\mR^d$, there exists a unique weak solution
to SDE \eqref{SDE} in the sense of Definition \ref{Def99} 
so that for each $T, R>0$ and $s,t\in[0,T]$,
\begin{align}\label{Def88}
\bE|A^b_{t\wedge\eta_R}-A^b_{s\wedge\eta_R}|^4\leq 
C_{T,R}|t-s|^{2(2-\alpha-\frac{d}{p})},
\end{align}
where  $\eta_R:=\inf\{t>0:|X_t|>R\}$. Moreover, $\bP\circ X^{-1}\in\sK^\alpha_p$ and
the conclusions in Theorem \ref{TH0} still hold.
\et 
\br\rm
In \cite[Theorem 2.6]{Ba-Ch1}, Bass and Chen require \eqref{Def88} to hold uniformly for $b_n$.
\er

\subsection{SDE with dissipative drift}
In this subsection, we consider SDE with dissipative drift but without  distributional part. 
First of all, we recall a stochastic Gronwall's inequality due to Scheutzew \cite{Sc} (see also \cite[Lemma 3.8]{Xi-Zh}).
\begin{lemma}[Stochastic Gronwall's inequality]\label{im}
Let $\xi(t)$ and $\eta(t)$ be two nonnegative c\`adl\`ag $\sF_t$-adapted processes, 
$A_t$ a continuous nondecreasing $\sF_t$-adapted process with $A_0=0$, $M_t$ a  local martingale with $M_0=0$. Suppose that
\begin{align}\label{Gron}
\xi(t)\leq\eta(t)+\int^t_0\xi(s)\dif A_s+M_t,\ \forall t\geq 0.
\end{align}
Then for any $0<q<p<1$ and stopping time $\tau$, we have
\begin{align}
\big[\mE(\xi(\tau)^*)^{q}\big]^{1/q}\leq \Big(\tfrac{p}{p-q}\Big)^{1/q}\Big(\mE \e^{pA_\tau/(1-p)}\Big)^{(1-p)/p}\mE\big(\eta(\tau)^*\big),  \label{gron}
\end{align}
where $\xi(t)^*:=\sup_{s\in[0,t]}\xi(s)$.
\end{lemma}
We first show the following result.
\bt\label{Le63}
Suppose that $b=b^{(1)}$ satisfies \eqref{Diss}, $\sigma$ satisfies \eqref{DG1} and for some $\gamma\in(0,1)$
\begin{align}\label{DG10}
|\sigma(x)-\sigma(y)|\leq c_1|x-y|^\gamma,\ \forall x,y\in\mR^d.
\end{align}
For each $x\in\mR^d$, there is one and only one element in $\sM^{0,\infty}_{\sigma, b}(x)$ denoted by $\mP_x$. 
Moreover, letting $\mE_x:=\mE^{\mP_x}$ and $P_t\varphi(x):=\mE_x\varphi(w_t)$, we have the following conclusions:
\begin{enumerate}[{\bf (a)}]
\item For each $m\in\mN$, there is a constant $C>0$ such that for all $T>0$ and $x\in\mR^d$,
\begin{align}\label{SD1}
\mE_x\left(\sup_{t\in[0,T]}|w_t|^m\right)+\mE_x\left(\int_0^T|w_s|^{m-1+\vartheta}\dif s\right)\leq C(1+|x|^m+T^m),
\end{align}
and for each $T>0$, there is a constant $C_{T}>0$ such that for all  $0\leq t_0<t_1\leq T$,
\begin{align}\label{Moment}
\mE_x |w_{t_1}-w_{t_0}|^{2m}\leq C_{T} |t_1-t_0|^m.
\end{align}
\item If $\vartheta=0$, then we have $P_t\varphi(x)=\int_{\mR^d} p(t,x,y)\varphi(y)\dif y$ with $p(t,x,y)$ satisfying \eqref{Es5} and \eqref{Es6}.
\item If $\vartheta>0$, then $P_t$ has a unique invariant probability measure.

\end{enumerate}
\et
\begin{proof}
{\bf (a)} For $n\in\mN$, let $\chi_n$ be the cutoff function defined by \eqref{Cut} and define
$$
b^n(x):=b(x)\cdot\chi_n(x).
$$
By \cite[Theorem 7.2.1]{St-Va}, there is a unique element $\mP^n_x\in \sM^{0,\infty}_{\sigma,b^n}(x)$.
Let $h(y):=(1+|y|^2)^{m/2}$. By \eqref{DG1} and
\eqref{Diss}, there are constants $\kappa_0',\kappa_1'>0$ independent of $n$ such that for all $|y|\leq n$,
$$
\sL^{\sigma} h(y)+\<\nabla h(y), b^n(y)\>\leq -\kappa'_0|y|^{m-1+\vartheta}+\kappa'_1(|y|^{m-1}+1).
$$
Let $\tau_n$ be defined by \eqref{Tau} with $R=n$. By definition, one has
\begin{align}\label{EG1}
\begin{split}
h(w_{t\wedge\tau_n})&=h(x)+\int_0^{t\wedge\tau_n}\Big(\sL^{\sigma} h+b^n\cdot\nabla h\Big)(w_s)\dif s+M_{t\wedge\tau_n}\\
&\leq h(x)-\kappa'_0\int_0^{t\wedge\tau_n}|w_s|^{m-1+\vartheta}\dif s+\kappa_1'\int^t_0(|w_{s\wedge\tau_n}|^{m-1}+1)\dif s+M_{t\wedge\tau_n},
\end{split}
\end{align}
where $M_{\cdot\wedge\tau_n}$ is a continuous martingale under $\mP^n_x$. 
By Lemma \ref{im}, for any $\delta\in(0,1)$, there is a constant $C>0$ such that for all  $T>0$ and $x\in\mR^d$,
\begin{align*}
\mE^{\mP^n_x}\left(\sup_{t\in[0,T]}h(w_{t\wedge\tau_n})^\delta\right)\leq Ch(x)^\delta
+C\left(\mE^{\mP^n_x}\int^T_0(|w_{s\wedge\tau_n}|^{m-1}+1)\dif s\right)^\delta.
\end{align*}
In particular, taking $\delta=1-\frac{1}{m}$ and by Young's inequality, we get
\begin{align*}
\mE^{\mP^n_x}\left(\sup_{t\in[0,T]}|w_{t\wedge\tau_n}|^{m-1}\right)
&\leq C(1+|x|^{m-1})
+CT^{1-\frac{1}{m}}\left(\mE^{\mP^n_x}\left(\sup_{t\in[0,T]}|w_{t\wedge\tau_n}|^{m-1}+1\right)\right)^{1-\frac{1}{m}}\\
&\leq C(1+|x|^{m-1})
+\frac{1}{2}\mE^{\mP^n_x}\left(\sup_{t\in[0,T]}|w_{t\wedge\tau_n}|^{m-1}+1\right)+CT^{m-1}.
\end{align*}
Here and below, the constant $C>0$ is independent of $n, T>0$ and $x\in\mR^d$.
Hence, for any $m\in\mN$,
$$
\mE^{\mP^n_x}\left(\sup_{t\in[0,T]}|w_{t\wedge\tau_n}|^{m-1}\right)\leq C(1+|x|^{m-1}+T^{m-1}),
$$
which implies
\begin{align}\label{EG9}
\mP^n_x(\tau_n\leq T)=\mP^n_x\left(\sup_{t\in[0,T]}|w_{t\wedge\tau_n}|>n\right)\leq C(1+|x|^m+T^m)/n^m\stackrel{n\to\infty}{\rightarrow} 0.
\end{align}
Thus, by \cite[p.250, Corollary 10.1.2]{St-Va}, there is a unique $\mP_x\in\sM^{0,\infty}_{\sigma,b}(x)$ such that for all $n\in\mN$,
\begin{align}\label{EG8}
\mP_x=\mP^n_x\ \mbox{ on}\ \ {\cB_{\tau_{n}}(\mC)}, 
\end{align}
and so,
\begin{align}\label{EG2}
\mE_x\left(\sup_{t\in[0,T]}|w_{t\wedge\tau_n}|^{m}\right)\leq C(1+|x|^m+T^m).
\end{align}
Substituting \eqref{EG2} into \eqref{EG1}, we also have
\begin{align}\label{EG3}
\mE_x\left(\int_0^{T\wedge\tau_n}|w_s|^{m-1+\vartheta}\dif s\right)\leq C(1+|x|^{m}+T^{m}).
\end{align}
By \eqref{EG9} and taking limits $n\to\infty$ for \eqref{EG2} and \eqref{EG3},  we get \eqref{SD1}.

On the other hand, notice that
$$
M_t:=w_t-y-\int^t_0b^n(w_s)\dif s\mbox{ is a continuous local martingale}
$$
with quadratic variation process 
$$
[M^i, M^j]_t:=\int^t_0 a^{jj}(w_s)\dif s.
$$
Hence, 
for each $T>0$ and $m\in\mN$, by Burkholder's inequality, there is a constant $C>0$ independent of $n$ such that for  all $0\leq t_0<t_1\leq T$,
\begin{align*}
\mE_x |w_{t_1}-w_{t_0}|^{2m}&\leq C\mE_x \left| \int_{t_0}^{t_1}b^n(w_s)\dif s\right|^{2m} +C \mE_x \left|M_{t_1}-M_{t_0}\right|^{2m}\\
&\leq C |t_1-t_0|^{2m-1}  \mE_x\int_{t_0}^{t_1} |b^n(w_s)|^{2m}\dif s
+ C\mE_x \left(\int_{t_0}^{t_1} |a(w_s)|\dif s\right)^m\\
&\leq C |t_1-t_0|^{2m-1}  
\mE_x\int_{t_0}^{t_1} (1+|w_s|^{2m\vartheta})\dif s
+C (t_1-t_0)^m\leq C |t_1-t_0|^{m}.
\end{align*}
\\
{\bf (b)} If $\vartheta=0$ in \eqref{Diss},  then $b$ is bounded measurable. Since $\sigma$ is uniformly non-degenerate and H\"older continuous,
it is well known that the semigroup $P_t$ admits a density $p(t,x,y)$
so that $P_t\varphi(x)=\int_{\mR^d}\varphi(y)p(t,x,y)\dif y$, and $p(t,x,y)$ satisfies \eqref{Es5} and \eqref{Es6}
(for example, see \cite[Theorem 1.1 and Section 4.2]{Ch-Hu-Xi-Zh}).
\\
\\
{\bf (c)} If $\vartheta>0$,  then by \eqref{SD1} with $m=1$ and the classical Bogoliubov-Krylov's argument, 
there exists an invariant probability measure associated with $P_t$. 
To show the uniqueness, as usual we show the strong Feller property and irreducibility of $P_t$. 
Let $\varphi$ be a bounded measurable function. For any $t>0$ and $x,y\in\mR^d$, we have
\begin{align}
\left|\mE_x\varphi(w_t)-\mE_y\varphi(w_t)\right|
&=\left|\mE_x (\varphi(w_t)1_{\{t<\tau_n\}})-\mE_y(\varphi(w_t)1_{\{t<\tau_n\}})\right|
+\|\varphi\|_\infty\Big(\mP_x(\tau_n\leq t)+\mP_y(\tau_n\leq t)\Big)\no\\
&\leq\left|\mE^{\mP^n_x}\varphi(w_t)-\mE^{\mP^n_y}\varphi(w_t)\right|+2\|\varphi\|_\infty\left(\mP^n_x(\tau_n\leq t)
+\mP^n_y(\tau_n\leq t)\right). \label{EC2}
\end{align}
By \eqref{EG9}, we have for any $R>0$,
$$
\lim_{n\to\infty}\sup_{|x|\leq R}\mP^n_x(\tau_n\leq t)=0.
$$
Since for each $t>0$, the mapping $x\mapsto \mE^{\mP^n_x}f(w_t)$ is continuous (see part {\bf (b)}),
letting $x\to y$ in \eqref{EC2} and then $n\to\infty$, we obtain the continuity of $x\mapsto P_tf(x)=\mE_x f(w_t)$.

Next we show the irreducibility of $P_t$. Fix $x_0, x\in \mR^d$, $t>0$ and $a>0$.
Choose $n$ large enough so that $D_n:=\{y: |y|<n\}$ contains $x_0$ and $B_a(x_0):=\{y: |y-x_0|\leq a\}$. Notice that
two-sided estimate \eqref{Es5} for the heat kernel of $P^n_t$ holds (see part {\bf (b)}). 
By \cite[Theorem 7.11]{Xi-Zh} or as in \cite[Theorem 2.2.4]{Ch-Zh}, one has
$$
\mP_{x} (w_t\in B_a(x_0))\geq\mP_{x}\big(w_{t}\in B_a(x_0); t<\tau_n)=\mP^n_{x}\big(w_{t}\in B_a(x_0); t<\tau_n)>0, 
$$
which means that $P_t$ is irreducible.
The proof is complete.
\end{proof}
Furthermore, we can show the following result.
\bl\label{Le510}
Assume {\bf (H$^\sigma_{\beta,q}$)} for some $\beta\in(0,1)$ and $q>d/\beta$, and $b=b^{(1)}$ satisfies \eqref{Diss}.
Let $\gamma\in(0,\beta]$ and $\nu>\frac{d}{1-\gamma}$.
For each $x\in\mR^d$, letting $\mP_x\in\sM^{0,\infty}_{\sigma,b}(x)$ be the unique element, 
there is a constant $C>0$ independent of $x$ such that for all $f\in C^\infty_c(\mR^d)$ and $T>0$,
\begin{align}\label{SD7}
\mE_x\left(\int^T_0f(w_s)\dif s\right)\leq C\|f\|_{-\gamma,\nu}\left(1+T+|x|\right).
\end{align}
Moreover, for any $T>0$ and $m\in\mN$, there is a constant $C_T>0$ such that for any $f\in C^\infty_c(\mR^d)$ and $0\leq t_0<t_1\leq T$,
\begin{align}\label{EJ1}
\mE_x\left|\int^{t_1}_{t_0}f(w_s)\dif s\right|^{2m}\leq C_{T}(t_1-t_0)^{(2-\gamma-\frac{d}{\nu})m}\|f\|^{2m}_{-\gamma,\nu}.
\end{align}
In particular, global Krylov's estimate holds for $\mP_x$, and $\sM^{0,\infty}_{\sigma,b}(x)\subset\sM^{\alpha,p}_{\sigma,b}(x)$.
\el
\begin{proof}
Under the assumptions of the theorem, by \eqref{RG1}, one sees that the assumptions of Theorem \ref{Le63} are satisfied.
For each $x\in\mR^d$, let $\mP_x\in\sM^{0,\infty}_{\sigma,b}(x)$ be the unique element. 
By Proposition \ref{Pr29}, there is a unique weak solution $(\Omega,\sF,(\sF_t)_{t\geq0},\bP; X,B)$ 
so that $\mP_x=\bP\circ X^{-1}$ and the following SDE is satisfied
\begin{align}\label{EG4}
X_t=x+\int^t_0\sigma(X_s)\dif B_s+\int^t_0b(X_s)\dif s.
\end{align}
Since $a=\sigma\sigma^*/2$ is H\"older continuous, by the classical Schauder theory of PDE (see \cite[p.56, Theorem 4.3.2]{Kr3}),
for any $f\in C^\infty_c(\mR^d)$,
there are $\lambda>0$ and a unique $u\in C^2_b(\mR^d)$ solving the following PDE:
$$
\sL^a u-\lambda u=f.
$$
Moreover, since the assumption (i) of Theorem \ref{Th44} is satisfied for $\alpha=\gamma$ and 
$p=\nu$, by \eqref{ER88} we also have 
\begin{align}\label{ER68}
\|u\|_{2-\gamma,\nu}\leq C\|f\|_{-\gamma,\nu}.
\end{align}
Now by \eqref{EG4} and It\^o's formula, we have
\begin{align}\label{SD2}
\begin{split}
u(X_t)&=u(x)+\int^{t}_0(\sL^{ a} u+b\cdot\nabla u)(X_s)\dif s+\int^{t}_0(\nabla u\cdot \sigma)(X_s)\dif B_s\\
&=u(x)+\int^{t}_0(f+\lambda u+b\cdot\nabla u)(X_s)\dif s
+\int^{t}_0(\nabla u\cdot\sigma)(X_s)\dif B_s.
\end{split}
\end{align}
Taking expectations and by \eqref{ER68}, Sobolev's embedding \eqref{Fr2}, we obtain
\begin{align*}
\bE\left(\int^{t}_0f(X_s)\dif s\right)&\leq (2+\lambda t)\|u\|_\infty+\|\nabla u\|_\infty\bE\left(\int^{t}_0| b(X_s)|\dif s\right)\\
&\leq C\|f\|_{-\gamma,\nu}\left(1+t+\bE\left(\int^{t}_0|X_s|^\vartheta\dif s\right)\right)
\stackrel{\eqref{SD1}}{\leq} C\|f\|_{-\gamma,\nu}\left(1+t+|x|\right).
\end{align*}
Thus, we get \eqref{SD7}.

On the other hand, for $0\leq t_0<t_1\leq T$, since
\begin{align}\label{GD2}
X_{t_1}-X_{t_0}=\int^{t_1}_{t_0} b(X_s)\dif s+\int^{t_1}_{t_0}\sigma(X_s)\dif B_s,
\end{align}
by \eqref{SD2} and easy calculations, we have
\begin{align}\label{GD7}
\begin{split}
\int^{t_1}_{t_0}f(X_s)\dif s
&=u(X_{t_1})-u(X_{t_0})-\int^{t_1}_{t_0}(\lambda u+ b\cdot\nabla u)(X_s)\dif s
-\int^{t_1}_{t_0}(\nabla u\cdot\sigma)(X_s)\dif B_s\\
&=(X_{t_1}-X_{t_0})\cdot\int^1_0\Big[\nabla u (rX_{t_1}+(1-r)X_{t_0})-\nabla u(X_{t_0})\Big]\dif r-\lambda_0\int^{t_1}_{t_0}u(X_s)\dif s\\
&-\int^{t_1}_{t_0} b(X_s)\cdot(\nabla u(X_s)-\nabla u(X_{t_0}))\dif s
-\int^{t_1}_{t_0}(\nabla u(X_s)-\nabla u(X_{t_0}))\cdot\sigma(X_s)\dif B_s.
\end{split}
\end{align}
Let $\delta:=1-\gamma-\frac{d}{\nu}$. By Sobolev's embedding \eqref{Fr2} and \eqref{ER68},
\begin{align}\label{Sob1}
|\nabla u(x)-\nabla u(y)|\lesssim \|\nabla u\|_{1-\gamma,\nu}|x-y|^\delta\lesssim\|f\|_{-\gamma,\nu}|x-y|^\delta.
\end{align}
By \eqref{GD7}, \eqref{Moment}, \eqref{Sob1},  
Burkholder's inequality 
and \eqref{ER68}, we have
\begin{align*}
\bE\left|\int^{t_1}_{t_0}f(X_s)\dif s\right|^{2m}
&\lesssim\|f\|^{2m}_{-\gamma,\nu}\bE|X_{t_1}-X_{t_0}|^{2m(1+\delta)}+\|f\|^{2m}_{-\gamma,\nu} |t_1-t_0|^{2m}\\
&\quad+\bE\left(\int^{t_1}_{t_0}|\nabla u(X_s)-\nabla u(X_{t_0})|^2\dif s\right)^{m}\\
&\lesssim\|f\|^{2m}_{-\gamma,\nu}| t_1-t_0| ^{m(1+\delta)}
+\|f\|^{2m}_{-\gamma,\nu}\bE\left(\int^{t_1}_{t_0}|X_{s}-X_{t_0}|^{2\delta}\dif s\right)^{m}\\ 
&\lesssim\|f\|^{2m}_{-\gamma,\nu}|t_1-t_0|^{m(1+\delta)}.
\end{align*}
Thus we complete the proof.
\end{proof}
\subsection{Proofs of Theorems \ref{TH0} and \ref{TH29}}

Let $\alpha\in(0,\frac{1}{2}]$, $p\in(\frac{d}{1-\alpha},\infty)$ and $\beta\in[\alpha,1]$, $q\in(\frac{d}{\beta},\infty)$.
Below we assume {\bf (H$^\sigma_{\beta,q}$)} and {\bf (H$^b_{\alpha,p}$)}. 
By (ii) of Theorem \ref{Th44}, there exists a constant $\lambda_0>0$ such that for all
$\lambda\geq \lambda_0$, there is a unique $\bu=\bu_\lambda: \mR^d\to\mR^d$ belonging to $H^{2-\alpha,p}$ so that
$$
(\sL^a-\lambda+b^{(2)}\cdot\nabla) \bu=-b^{(2)}\ \mbox{in }\ H^{-\alpha,p}.
$$
By \eqref{ER88}, for any $\theta\in[0,2]$, 
there is a constant $C>0$ such that for all $\lambda\geq\lambda_0$,
\begin{align}\label{ER98}
\lambda^{1-\frac{\theta}{2}} \|\bu\|_{\theta-\alpha,p}\leq C\|b^{(2)}\|_{-\alpha,p}.
\end{align}
In particular, taking $\theta\in(1+\alpha+\frac{d}{p},2)$ and by Sobolev's embedding \eqref{Fr2}, we can choose $\lambda$ large enough so that
\begin{align}\label{ER67}
\|\nabla \bu\|_\infty\leq 1/2.
\end{align}
Now, define 
$$
\Phi(x):=x+\bu(x): \mR^d\to\mR^d.
$$
By \eqref{ER67} and \eqref{ER98} with $\theta=2$, it is easy to see that 
\begin{align}\label{FG1}
\tfrac{1}{2}|x-y|\leq |\Phi(x)-\Phi(y)|\leq 2|x-y|,\ \  \|\mI-\nabla\Phi\|_{1-\alpha,p}=\|\nabla \bu\|_{1-\alpha,p}\leq C\|b^{(2)}\|_{-\a,p}.
\end{align}
Hence,  $\Phi\in\cD^{1-\alpha}_p$  and
\begin{align}\label{Nee}
\sL^\sigma\Phi+b^{(2)}\cdot\nabla\Phi=\lambda \bu\ \ \mbox{ in $H^{-\alpha,p}$}.
\end{align}
Define
\begin{align}\label{Ne}
\tilde\sigma:=(\nabla\Phi\cdot\sigma)\circ\Phi^{-1},\quad
 \tilde b:=(\lambda \bu+b^{(1)}\cdot\nabla\Phi)\circ\Phi^{-1}.
\end{align}
We have the following key observation.
\bl\label{Le52}
For $\lambda$ large enough, there are $\tilde\kappa_0,\tilde\kappa_1,\tilde\kappa_2>0$ such that for all $y\in\mR^d$,
\begin{align}\label{Disp}
\frac{\<y, \tilde b(y)\>}{\sqrt{1+|y|^2}}\leq -\tilde\kappa_0|y|^{\vartheta}+\tilde\kappa_1
\quad\text{and}\quad|\tilde b(y)|\leq \tilde\kappa_2(1+|y|^{\vartheta}),
\end{align}
where $\vartheta$ is the same as in \eqref{Diss}.
Moreover, $\tilde\sigma$ satisfies {\bf (H$^\sigma_{\beta', q'}$)} with 
$\beta'=\beta\wedge (1-\alpha)$ and $q'$ being defined by \eqref{GS1}.
\el
\begin{proof}
If $\vartheta=0$, there is nothing to prove \eqref{Disp}. Below we assume $\vartheta>0$.
First of all, it is clear that
$$
|\tilde b(y)|\leq\lambda\|\bu\|_\infty+\kappa_2(1+|\Phi^{-1}(y)|^\vartheta)\|\nabla\Phi\|_\infty\leq\tilde\kappa_2(1+|y|^\vartheta).
$$
Observing that
$$
y=\Phi^{-1}(y)+ \bu\big(\Phi^{-1}(y)\big),\ \ \nabla\Phi(x)=\mI+\nabla\bu(x),
$$
by the definition of $\tilde b$ and \eqref{Diss},
we have
\begin{align*}
\frac{\<y, \tilde b(y)\>}{\sqrt{1+|y|^2}}&=\frac{\lambda\<y, \bu(\Phi^{-1}(y))\>}
{\sqrt{1+|y|^2}}+\frac{\<y,b^{(1)}(\Phi^{-1}(y))\>}{\sqrt{1+|y|^2}}+\frac{\<y,(b^{(1)}\cdot\nabla \bu)(\Phi^{-1}(y))\>}{\sqrt{1+|y|^2}}\\
&\leq \lambda\| \bu\|_{\infty} +\frac{\big\langle\Phi^{-1}(y),  b^{(1)}\big(\Phi^{-1}(y)\big)\big\rangle}{\sqrt{1+|y|^2}}
+\frac{|b^{(1)}\big(\Phi^{-1}(y)\big)|(\|\bu\|_{\infty}+\|\nabla \bu\|_\infty|y|)}{\sqrt{1+|y|^2}}\\
&\leq \lambda\|\bu\|_\infty+\left(\kappa_2-\kappa_1|\Phi^{-1}(y)|^{\vartheta}\right)\frac{\sqrt{1+|\Phi^{-1}(y)|^2}}{\sqrt{1+|y|^2}}
+\kappa_3(1+|\Phi^{-1}(y)|^\vartheta)\frac{\|\bu\|_{\infty}+\|\nabla\bu\|_\infty|y|}{\sqrt{1+|y|^2}}.
\end{align*}
By \eqref{ER98} and Sobolev's embedding \eqref{Fr2}, we have $\lim_{\lambda\to\infty}\|\nabla\bu_\lambda\|_\infty=0$.
The first estimate in \eqref{Disp} follows by choosing $\lambda$ large enough and \eqref{FG1}.
Moreover, by (i) of Proposition \ref{Pr24}, $\tilde\sigma$ satisfies {\bf (H$^\sigma_{\beta', q'}$)}.
\end{proof}

Now we can give 
\begin{proof}[Proof of Theorem \ref{TH0}]
For each $x\in\mR^d$,  by Proposition \ref{Pr24}, Lemmas \ref{Le52}, \ref{Le510} and Theorem \ref{Le63}, the unique martingale solution 
$\mP_x\in\sM^{\alpha,p}_{\sigma,b}(x)$ is given by
\begin{align}\label{Ma}
\mP_x=\tilde\mP_{\Phi(x)}\circ\Phi,
\end{align}
where $\tilde\mP_{y}\in\sM^{\alpha,p}_{\tilde\sigma,\tilde b}(y)$ is the unique martingale solution
starting from $y$ associated with $\tilde\sigma,\tilde b$. We shall write $\tilde\mE_y:=\mE^{\tilde\mP_y}$
and $\tilde P_t\varphi(y):=\tilde\mE_y\varphi(w_t)$.

(i) It follows by \eqref{Ma} and Theorem \ref{Le63}, Lemma \ref{Le510}.

(ii) If $\vartheta=0$ in \eqref{Diss},  then by \eqref{Ma}, Lemma \ref{Le52} and (ii) of Theorem \ref{Le63}, we have
$$
p(t,x,y)=\tilde p(t,\Phi(x),\Phi(y))\det(\nabla\Phi(y)),
$$
where $\tilde p(t,x,y)$ is the heat kernel of $\tilde P_t$.
Since $\tilde p(t,x,y)$ enjoys the estimates \eqref{Es5} and \eqref{Es6}, by \eqref{FG1}, it is easy to see that
$p(t,x,y)$ also enjoys the estimates \eqref{Es5} and \eqref{Es6}.

Next we show the probabilistic representation part.
For any $\varphi\in \cap_{m\in\mN}H^{m,p}$, by Theorem \ref{Th404} and Remark \ref{Re47},
there is a unique $u^\varphi\in \mH^{2-\alpha,p}_T\cap\mH^{2,p}_T$ satisfies
\begin{align}\label{EC1}
u^\varphi(t)=\varphi+\int^t_0(\sL^{\tilde a}+\tilde b\cdot\nabla)u^\varphi(s)\dif s.
\end{align}
Since $u^\varphi\in \mH^{2,p}_T$, by applying generalized It\^o's formula to $(t,x)\mapsto u^\varphi(T-t,x)$ (see \cite[p.121]{Kr-1}), 
we get $\tilde\mE_x  u^\varphi(0, w_T)=u^\varphi(T,x)$, i.e. $\tilde\mE_x  \varphi(w_T)=u^\varphi(T,x).$ 
Furthermore, for general $\varphi\in H^{2-\alpha,p}$, let $\varphi_n:=\varphi*\rho_n$ be the mollifying approximation.
By Theorem \ref{Th404}, it is easy to see that
$$
\|\p_tu^{\varphi_n}-\p_tu^\varphi\|_{\mH^{-\alpha,p}_T}+\|u^{\varphi_n}-u^\varphi\|_{\mH^{2-\alpha,p}_T}
\leq C\|\varphi_n-\varphi\|_{2-\alpha,p}\stackrel{n\to\infty}{\to} 0.
$$
By taking limits for $u^{\varphi_n}(t,x)=\tilde\mE_x  \varphi_n(w_t)$, 
we get the probabilistic representation  for the unique solution $u^\varphi$ of \eqref{EC1}:
\begin{align}\label{Prob}
u^\varphi(t,x)=\tilde\mE_x  \varphi(w_t)\in \mH^{2-\alpha,p}_T,\ \varphi\in H^{2-\alpha,p}.
\end{align}
Moreover, by \eqref{Ma} we have
$$
P_t\varphi (x)=\mE^{\mP_x}\varphi (w_t)=\mE^{\tilde\mP_{\Phi(x)}\circ\Phi}\varphi (w_t)=(\tilde P_t (\varphi \circ\Phi^{-1}))(\Phi(x)),
$$
Since $\Phi\in\cD^{1-\alpha}_p$,  by  Lemmas \ref{Prod}, \eqref{Le211} and Proposition \ref{Pr22}, we have
\begin{align*}
\|\varphi \circ \Phi^{-1}\|_{2-\a,p}\lesssim&\|\varphi \circ \Phi^{-1}\|_{p}+ \|\nabla (\varphi \circ \Phi^{-1})\|_{1-\a, p}
=\|\varphi \|_p+  \|\nabla \varphi \circ \Phi^{-1}\cdot \nabla \Phi^{-1} \|_{1-\a,p}\\
\lesssim& \|\varphi \|_p+ \|\nabla \varphi \|_{1-\a,p}( \|\nabla \Phi^{-1}-\mI\|_{1-\a, p}+1) \lesssim \|\varphi \|_{2-\a, p}.
\end{align*}
Hence,  by \eqref{Prob}, $\tilde P_t(\varphi \circ \Phi^{-1})\in\mH^{2-\a, p}_T$ and
$P_t \varphi=\tilde P_t(\varphi \circ \Phi^{-1})\circ\Phi\in\mH^{2-\a, p}_T$. 
By \eqref{YP1}, \eqref{YP2} and \eqref{EC1}, one sees that $P_t\varphi$ satisfies \eqref{PDE9}. 

(iii) If $\vartheta>0$, by Theorem \ref{Le63}, $\tilde P_t$ admits a unique invariant probability measure $\tilde\mu$,
and by Lemma \ref{Le510}, there is a constant $C$ independent of $T$ and $y$ such that for 
any $\gamma\in(0,\beta']$ 
, $\nu>\frac{d}{1-\gamma}$ 
and any $f\in C_c^\infty(\mR^d)$, $y\in\mR^d$,
$$
\frac{1}{T} \int_0^T \tilde P_tf(y)\dif t=\frac{1}{T}\tilde\mE_y\left(\int^T_0f(w_s)\dif s\right)\leq C\|f\|_{-\gamma,\nu}\frac{1+|y|+T}{T},
$$
which implies by Birkhorff's ergodicity theorem,
$$
\tilde\mu(f)=\lim_{T\to\infty} \frac{1}{T} \int_0^T \tilde P_t f(y)\dif t\leq C\|f\|_{-\gamma,\nu}.
$$
Hence, $\tilde\mu(\dif y)=\tilde\varrho(y)\dif y$ and
$\tilde\varrho\in H^{\gamma,r}$, where $r=\tfrac{\nu}{\nu-1}\in (1, \tfrac{d}{d+\gamma-1})$.
Finally, by \cite[Proposition 2.8]{Xi-Zh}, $\mu:=\tilde\mu\circ\Phi$ is the unique invariant probability measure of
$P_t$. Moreover, it is clear that $\mu(\dif x)=\varrho(x)\dif x$ with $\varrho(x)=\tilde\varrho\circ\Phi(x)\det(\nabla\Phi(x))$.
Noticing that  $\tfrac{1}{p}-\tfrac{1-\a-\gamma}{d} < \tfrac{d+\gamma-1}{d}<\tfrac{1}{r}$, 
by Sobolev's embedding, Lemmas \ref{Prod} and \ref{Le211}, one sees that $\varrho\in H^{\gamma,r}$.
\end{proof}

Finally we prove Theorem \ref{TH29}.

\begin{proof}[Proof of Theorem \ref{TH29}]
By Proposition \ref{Pr29} and Theorem \ref{TH0}, it suffices to show that  for any weak solution $(\Omega,\sF, (\sF_t)_{t\geq 0}, \bP; X,B)$, if \eqref{Def88} holds, then the local Krylov estimate holds for the law of $X$.
Let $f\in C^\infty_c(\mR^d)$. By Theorem \ref{Th44}, for $\lambda$ large enough, there is a unique $u\in H^{2-\alpha,p}$ solving the following PDE:
$$
\sL^au -\lambda u=f.
$$
By \eqref{Fr2} and \eqref{ER88}, we have
\begin{align}\label{UQ1}
|\nabla u(x)-\nabla u(y)|\lesssim\|\nabla u\|_{1-\alpha,p}|x-y|^{1-\a-d/p}\lesssim \|f\|_{-\alpha,p}|x-y|^{1-\a-d/p}.
\end{align}
Moreover, by Schauder's theory of PDE (see \cite[p.56, Theorem 4.3.2]{Kr3}), we also have $u\in C^2$. Hence,
by It\^o's formula for Dirichlet processes (see Lemma \ref{Ito}),
$$
u(X_t)=u(x)+\int^t_0
(\lambda u-f)(X_s)\dif s
+\int^t_0\nabla u(X_s)\cdot\dif A^b_s+\int^t_0(\nabla u\cdot\sigma)(X_s)\dif B_s.
$$
As in the calculations of \eqref{GD7}, for $t_0<t_1$, we have
\begin{align}\label{GQ3}
\begin{split}
\int^{t_1}_{t_0}f(X_s)\dif s
&=(X_{t_1}-X_{t_0})\int^1_0\Big[\nabla u(rX_{t_1}+(1-r)X_{t_0})-\nabla u(X_{t_0})\Big]\dif r-\lambda\int^{t_1}_{t_0} u(X_s)\dif s\\
&+\int^{t_1}_{t_0}(\nabla u(X_s)-\nabla u(X_{t_0}))\cdot\dif A^b_s
-\int^{t_1}_{t_0}(\nabla u(X_s)-\nabla u(X_{t_0}))\cdot\sigma(X_s)\dif B_s.
\end{split}
\end{align}
Since $\bE|A^b_{t_1\wedge\eta_R}-A^b_{t_0\wedge\eta_R}|^4\leq C|t_1-t_0|^{2(2-\alpha-\frac{d}{p})}$, by \eqref{Def9} we have
$$
\bE|X_{t_1\wedge\eta_R}-X_{t_0\wedge\eta_R}|^4\leq C|t_1-t_0|^{2},
$$
and by \eqref{UQ1},
$$
\bE|\nabla u(X_{t_1\wedge\eta_R})-\nabla u(X_{t_0\wedge\eta_R})|^4 \leq C\|f\|_{-\alpha,p}^4|t_1-t_0|^{2(1-\a-d/p)}.
$$
Thus, by $\frac{1}{2}-\alpha-\frac{d}{p}>0$ and Lemma \ref{Le72} with $p=q=4$, we get 
$$
\bE\left|\int^{t_1\wedge \eta_R}_{t_0\wedge\eta_R}(\nabla u(X_s)-\nabla u(X_{t_0}))\cdot\dif A^b_s\right|^2\leq C\|f\|_{-\alpha,p}^2|t_1-t_0|^{2-\alpha-\frac{d}{p}}. 
$$
By \eqref{GQ3} and as in proving  \eqref{EJ1}, we have
\begin{align*}
\bE\left|\int^{t_1\wedge\eta_R}_{t_0\wedge\eta_R}f(X_s)\dif s\right|^2&\leq C\|f\|_{-\alpha,p}^2|t_1-t_0|^{2-\alpha-\frac{d}{p}}.
\end{align*}
The proof is complete.
\end{proof}


\begin{thebibliography}{999}

\bibitem{Ba-Ch-Da}Bahouri H., Chemin J.Y. and Danchin R.: {\it Fourier Analysis and Nonlinear Partial Differential Equations}.
Springer-Verlag, 2011.

\bibitem{Ba-Ch1}Bass R.F. and Chens Z.Q.: Stochastic differential equations for Dirichlet processes. 
{\it Probability Theory and Related Fields}. {\bf 121}(3) (2001), pp.422-446. 

\bibitem{Ba-Ch2}Bass R.F. and Chen Z.Q.: Brownian motion with singular drift. {\it The Annals of Probability}. {\bf31}(2) (2003), 791-817.	

\bibitem{Bo-Kr-Ro-Sa}Bogachev V.I., Krylov N.V., R\"ockner M., Shaposhnikov S.V.:
{\it Fokker-Planck-Kolmogorov Equations},\ AMS Mathematical Surveys and Monographs 207,
Publisher: American Mathematical Society, Year: 2015.

\bibitem{Ch-Hu-Xi-Zh}Chen Z.Q., Hu E., Xie L. and Zhang X.: Heat kernels for non-symmetric diffusion operators with jumps.
{\it J. Differential Equations}. {\bf 263} (2017) 6576-6634.

\bibitem{Ch-Zh0}Chen Z.Q. and Zhang X.: Uniqueness of stable-like processes. arXiv:1604.02681.

\bibitem{Ch-En}Cherny A.S. and Engelbert H.J.: Singular Stochastic Differential Equations. {\it Lecture Notes in Mathematics}, {\bf 1858}, Springer, 2005.

\bibitem{Ch-Zh}Chung K.L. and Zhao Z.: {\it From Brownian motion to Schr\"odinger’s equation}. 
 Springer-Verlag, Berlin, 1995. 




\bibitem{Fl-Is-Ru}Flandoli F., Issoglio E. and Russo F.: Multidimensional stochastic differential equations with distributional drift. 
{\it Transactions of the American Mathematical Society}. {\bf 369}(3), (2017), 1665-1688.

\bibitem{Fl-Ru-Wo1} Flandoli F., Russo F. and Wolf J.: Some SDEs with distributional drift. I. General calculus. {\it Osaka J. Math.}, {\bf 40}(2):493-542, 2003.

\bibitem{Fl-Ru-Wo2} Flandoli F., Russo F. and Wolf J.: Some SDEs with distributional drift. II. Lyons- Zheng structure, 
It\^o's formula and semimartingale characterization. {\it Random Oper. Stochastic Equations}, {\bf 12}(2):145-184, 2004.

\bibitem{Fo}F\"ollmer H.: Calcul d'It\^o sans probabilit\'es, {\it S\'eminaire de Probabilit\'es XV}, 143-150,  Springer, Berlin (1981). 
	
\bibitem{Gu-Im-Pe}Gubinelli M., Imkeller P. and Perkowski N.: Paracontrolled distributions and singular PDEs. {\it Forum of Mathematics}, (2015).
	
\bibitem{He-Zh}He K. and Zhang X:  One dimensional stochastic differential equations with distributional drifts. 
{\it Acta Mathematicae Applicatae Sinica}, English Series. {\bf 23}(3), (2007), 501-512. 	

\bibitem{Hu-Le-Mi}Hu Y., L\^e K. and Mytnik L.: Stochastic differential equation for Brox diffusion. 
{\it Stochastic Processes and their Applications}, {\bf 127} (2017), 2281-2315.

\bibitem{Kr-1}Krylov N.V.:  {\it Controlled Diffusion Processes}.  Springer-Verlag, 1980. 

\bibitem{Kr1}Krylov N.V.:  A generalization of the Littlewood-Paley inequality and some other results related to stochastic partial differential equations. 
{\it Ulam Quart}. {\bf  2}(4), (1994), p.16.


\bibitem{Kr0}Krylov N.V.:   The heat equation in $L^q((0,T); L^p)$-spaces with weights. {\it SIAM J. Math. Anal.}, {\bf 32}, (2001), 1117-1141. 

\bibitem{Kr3} Krylov N.V.: {\it Lectures on elliptic and parabolic equations in H\"older spaces}. American Mathematical Soc. { \bf Vol. 12} (1996). 
 

\bibitem{Kr2} Krylov N.V.: {\it Lectures on elliptic and parabolic equations in Sobolev spaces}. American Mathematical Soc. {\bf Vol. 96} (2008). 




\bibitem{Kr-Ro}Krylov N.V. and R\"ockner M.:
Strong solutions of stochastic equations with singular time dependent drift.
{\it Probab. Theory Relat. Fields}, {\bf 131} (2005), 154-196.

\bibitem{Ru-Tr}Russo F. and Trutnau G.: Some parabolic PDEs whose drift is an irregular
random noise in space. {\it The Annals of Probability}, Vol. {\bf 35}, No. 6, (2007) 2213-2262.

\bibitem{Re-Yo}Revu D. and Yor M.: {\it Continuous martingale and Brownian motion}. Third Edition, Springer-Verlag, Berlin, 1999. 

\bibitem{Sc}Scheutzow M.: A stochastic Gronwall's lemma. {\it Infinite Dimensional Analysis, Quantum Probability 
and Related Topics}, {\bf 16}, No. 2 (2013) 1350019 (4 pages).

\bibitem{St-Va}Stroock D.W. and Varadhan S.S.: {\it Multidimensional diffusion processes}. Springer-Verlag, Berlin, 1979. 	


\bibitem{Xi-Zh}Xie L. and Zhang X.: Ergodicity of stochastic differential equations with jumps and singular coefficients. {\it arXiv:1705.07402}.

\bibitem{ZhQ}Zhang Q.S.: Gaussian bounds for the fundamental solutions of $\nabla(A\nabla u)+B\nabla u-u_t =0$. {\it Manuscripta Math.}
93(3) (1997) 381-390.

\bibitem{Zh0}Zhang X.: Stochastic homeomorphism flows of
SDEs with singular drifts and Sobolev diffusion coefficients.
{\it Electron. J. Probab. \bf 16} (2011), 1096-1116.

\bibitem{Zh1}Zhang X.: Stochastic differential equations with Sobolev diffusion and singular drift. {\it Annals
of Applied Probability}. {\bf 26}, No. 5, 2697-2732(2016).


\bibitem{Zv}Zvonkin, A.K.: A transformation of the phase space of a diffusion process
that removes the drift.  {\it Mat. Sbornik}. {\bf 93} (135) (1974), 129-149.



\end{thebibliography}
\end{document}